# Nonvanishing of quadratic Dirichlet $L$-functions at $s = \frac{1}{2}$

By K. Soundararajan*

## 1. Introduction

The Generalized Riemann Hypothesis (GRH) states that all nontrivial zeros of Dirichlet $L$-functions lie on the line $\text{Re}(s) = \frac{1}{2}$. Further, it is believed that there are no $\mathbb{Q}$-linear relations among the nonnegative ordinates of these zeros. In particular, it is expected that $L(\frac{1}{2}, \chi) \neq 0$ for all primitive characters $\chi$, but this remains unproved. It appears to have been conjectured first by S. D. Chowla [5] in the case when $\chi$ is a quadratic character. In addition to numerical evidence (see [16] and [17]) the philosophy of N. Katz and P. Sarnak [13] lends theoretical support to this belief. Assuming the GRH, they proved that (oral communication) for at least $(19 - \cot(\frac{1}{4}))/16 > \frac{15}{16}$ of the fundamental discriminants $|d| \leq X$, $L(\frac{1}{2}, (\frac{d}{\cdot})) \neq 0$. Independently, A. E. Özluk and C. Snyder [15] showed, also assuming GRH, that $L(\frac{1}{2}, \chi_d) \neq 0$ for at least $\frac{15}{16}$ of the fundamental discriminants $|d| \leq X$. Katz and Sarnak also developed conjectures on the low-lying zeros in this family of $L$-functions (analogous to the Pair Correlation conjecture regarding the vertical distribution of zeros of $\zeta(s)$) which imply that $L(\frac{1}{2}, (\frac{d}{\cdot})) \neq 0$ for almost all fundamental discriminants $d$. In a different vein, R. Balasubramanian and V. K. Murty [1] showed that for a (small) positive proportion of the characters $\pmod q$, $L(\frac{1}{2}, \chi) \neq 0$. Recently, H. Iwaniec and P. Sarnak [10] have demonstrated that this proportion is at least one third.

For integers $d \equiv 0$, or $1 \pmod 4$ we put $\chi_d(n) = (\frac{d}{n})$. Notice that $\chi_d$ is a real character with conductor $\leq |d|$. If $d$ is an odd, positive, square-free integer then $\chi_{8d}$ is a real, primitive character with conductor $8d$, and with $\chi_{8d}(-1) = 1$. In [19], we considered the family of quadratic twists of a fixed Dirichlet $L$-function $L(s, \psi)$. Precisely, we considered the family $L(s, \psi \otimes \chi_{8d})$ for odd, positive, square-free integers $d$. When $\psi$ is *not* quadratic we showed that at least $\frac{1}{5}$ of these $L$-functions are not zero at $s = \frac{1}{2}$, and indicated how this proportion may be improved to $\frac{1}{3}$. The most interesting case when

*Research supported in part by the American Institute of Mathematics (AIM).

$\psi = 1$ (or, what amounts to the same, when $\psi$ is quadratic) turned out to be substantially different from the case when $\psi$ was not quadratic. There arose here an "off-diagonal" contribution which we were unable to evaluate in [19]. In this paper, we resolve the case when $\psi = 1$ and establish that a high proportion of quadratic Dirichlet $L$-functions are not zero at $s = \frac{1}{2}$.

THEOREM 1. *For at least 87.5% of the odd square-free integers $d \geq 0$, $L(\frac{1}{2}, \chi_{8d}) \neq 0$. Precisely, for all large $x$, and any fixed $\varepsilon > 0$,*

$$\sum_{\substack{d \leq x \\ L(\frac{1}{2},\chi_{8d}) \neq 0}} \mu(2d)^2 \geq \left(\frac{7}{8} - \varepsilon\right) \sum_{d \leq x} \mu(2d)^2.$$

It is striking that the proportion of nonvanishing in Theorem 1 is more than twice as good as the proportion obtained when $\psi$ is not quadratic, and also the proportion obtained by Iwaniec and Sarnak in the family of all Dirichlet $L$-functions (mod $q$). One explanation for this is that if $L(\frac{1}{2}, \chi_{8d}) = 0$ then automatically $L'(\frac{1}{2}, \chi_{8d}) = 0$; this does not hold in the other two families. This makes it more unlikely for $L(s, \chi_{8d})$ to vanish at $\frac{1}{2}$ than in the other cases. Another explanation is provided by the Katz-Sarnak models [13]. The zeros of $L(s, \chi_{8d})$ are governed by a symplectic law where there is greater repulsion of $s = \frac{1}{2}$, whereas the zeros of the $L(s, \psi \otimes \chi_{8d})$ ($\psi$ not quadratic) and $L(s, \chi)$ ($\chi \pmod{q}$) are governed by a unitary law with no repulsion of $s = \frac{1}{2}$. The same proportion $\frac{7}{8}$ appears in work of E. Kowalski and P. Michel [14] concerning the rank of $J_0(q)$. They showed that the proportion of odd, primitive, modular forms $f$ of weight 2 and level $q$ with $L'(f, \frac{1}{2}) \neq 0$ is at least $\frac{7}{8}$ (note that since $f$ is odd, $L(f, \frac{1}{2}) = 0$). This coincidence may be 'explained' by noting that the Kowalski-Michel family is governed by an odd orthogonal symmetry (SO($2N + 1$)) and the distribution of the second eigenvalue in such a family matches precisely the distribution of the first eigenvalue in the symplectic family of Theorem 1 (see pages 10–15 of [13]).

In Theorem 1 we considered only fundamental discriminants divisible by 8. We may replace this by fundamental discriminants in any arithmetic progression $a \pmod{b}$; this would include all the quadratic twists of $\psi$ for any quadratic character $\psi$. Also, the point $\frac{1}{2}$ is not special. A similar result (with a different proportion) may be established for any point $\sigma + it$ in the critical strip.

Earlier work of Jutila [12] shows that that there are $\gg X/\log X$ fundamental discriminants $d$ with $|d| \leq X$ such that $L(\frac{1}{2}, \chi_d) \neq 0$. He achieved this by evaluating the first and second moments of $L(\frac{1}{2}, \chi_d)$. That is, for two positive constants $c_1$ and $c_2$ he established

$$\text{(1.1)} \qquad \sum_{|d| \leq X} L(\tfrac{1}{2}, \chi_d) \sim c_1 X \log X,$$

and

$$\text{(1.2)} \qquad \sum_{|d| \leq X} |L(\tfrac{1}{2}, \chi_d)|^2 \sim c_2 X (\log X)^3,$$

where $d$ ranges over fundamental discriminants in both sums above. By Cauchy's inequality it follows that the number of fundamental discriminants $|d| \leq X$ such that $L(\tfrac{1}{2}, \chi_d) \neq 0$ exceeds the ratio of the square of the quantity in (1.1) to the quantity in (1.2) which is $\gg X/\log X$.

The improvement in Theorem 1 comes from the introduction of a "mollifier." Historically mollifiers appear first in work of Bohr and Landau [2] on zeros of the Riemann zeta function. Later this idea was used with remarkable success by Selberg [18] to demonstrate that a positive proportion of the zeros of $\zeta(s)$ lie on the critical line. Our aim here is to find a mollifier

$$\text{(1.3)} \qquad M(d) = \sum_{l \leq M} \lambda(l) \sqrt{l}\left(\frac{8d}{l}\right),$$

such that the mollified first and second moments are comparable. Precisely, we want

$$\sum_{d \leq x} \mu(2d)^2 L(\tfrac{1}{2}, \chi_{8d}) M(d) \asymp \sum_{d \leq x} \mu(2d)^2 |L(\tfrac{1}{2}, \chi_{8d}) M(d)|^2 \asymp x.$$

By Cauchy's inequality this demonstrates that a positive proportion of odd square-free $d$'s satisfy $L(\tfrac{1}{2}, \chi_{8d}) \neq 0$. In Section 6 we achieve this by choosing an optimal mollifier which has the shape (for an odd integer $l \leq M$)

$$\lambda(l) \qquad \text{roughly proportional to} \qquad \frac{\mu(l)}{l} \frac{\log^2(M/l)}{\log^2 M} \frac{\log(X^{\frac{3}{2}} M^2 l)}{\log M}.$$

By taking $M = X^{\frac{1}{2}-\varepsilon}$ and evaluating the first and second mollified moments for this optimal choice, we prove Theorem 1.

We now give a detailed outline of the proof of Theorem 1. Let $\{f_n\}_{n=1}^{\infty}$ be any sequence of complex numbers and let $F$ denote a nonnegative Schwarz class function compactly supported in the interval $(1,2)$. We define

$$\mathcal{S}(f_d; F) = \mathcal{S}_X(f_d; F) = \frac{1}{X} \sum_{d \text{ odd}} \mu^2(d) f_d F\left(\frac{d}{X}\right).$$

Let $Y > 1$ be a real parameter to be chosen later and write $\mu^2(d) = M_Y(d) + R_Y(d)$ where

$$M_Y(d) = \sum_{\substack{l^2 \mid d \\ l \leq Y}} \mu(l), \quad \text{and} \quad R_Y(d) = \sum_{\substack{l^2 \mid d \\ l > Y}} \mu(l).$$

Define

$$\mathcal{S}_M(f_d; F) = \mathcal{S}_{M,X,Y}(f_d; F) = \frac{1}{X} \sum_{d \text{ odd}} M_Y(d) f_d F\left(\frac{d}{X}\right),$$

and

$$\mathcal{S}_R(f_d; F) = \mathcal{S}_{R,X,Y}(f_d; F) = \frac{1}{X} \sum_{d \text{ odd}} |R_Y(d) f_d| F\left(\frac{d}{X}\right),$$

so that $\mathcal{S}(f_d; F) = \mathcal{S}_M(f_d; F) + O(\mathcal{S}_R(f_d; F))$.

In this notation, we seek to evaluate the mollified moments $\mathcal{S}(L(\tfrac{1}{2}, \chi_{8d}) M(d); \Phi)$ and $\mathcal{S}(|L(\tfrac{1}{2}, \chi_{8d}) M(d)|^2; \Phi)$. Here, and in the sequel, $\Phi$ is a smooth Schwarz class function compactly supported in $(1, 2)$ and we assume that $0 \leq \Phi(t) \leq 1$ for all $t$. For integers $\nu \geq 0$ we define

$$\Phi_{(\nu)} = \max_{0 \leq j \leq \nu} \int_1^2 |\Phi^{(j)}(t)| dt.$$

For any complex number $w$ define

$$\check{\Phi}(w) = \int_0^\infty \Phi(y) y^w dy,$$

so that $\check{\Phi}(w)$ is a holomorphic function of $w$. Integrating by parts $\nu$ times, we get that

$$\check{\Phi}(w) = \frac{1}{(w+1)\cdots(w+\nu)} \int_0^\infty \Phi^{(\nu)}(y) y^{w+\nu} dy,$$

so that for $\text{Re}(w) > -1$ we have

(1.4) $$|\check{\Phi}(w)| \ll_\nu \frac{2^{\text{Re}(w)}}{|w+1|^\nu} \Phi_{(\nu)}.$$

To evaluate these moments, we first need "approximate functional equations" for $L(\tfrac{1}{2}, \chi_{8d})$ and $|L(\tfrac{1}{2}, \chi_{8d})|^2$. For integers $j \geq 1$ put $\omega_j(0) = 1$ and for $\xi > 0$ define $\omega_j(\xi)$ by

(1.5) $$\omega_j(\xi) = \frac{1}{2\pi i} \int_{(c)} \left(\frac{\Gamma(\tfrac{s}{2} + \tfrac{1}{4})}{\Gamma(\tfrac{1}{4})}\right)^j \xi^{-s} \frac{ds}{s}$$

where $c$ is any positive real number. Here, and henceforth, $\int_{(c)}$ stands for $\int_{c-i\infty}^{c+i\infty}$. In Lemma 2.1, we shall show that $\omega_j(\xi)$ is a real-valued smooth function on $[0, \infty)$ and that $\omega_j(\xi)$ decays exponentially as $\xi \to \infty$. As usual, $d_j(n)$ will denote the $j$-th divisor function; that is the coefficient of $n^{-s}$ in the Dirichlet series expansion of $\zeta(s)^j$. For integers $j \geq 1$, we define

$$(1.6) \qquad A_j(d) = \sum_{n=1}^{\infty} \left(\frac{8d}{n}\right) \frac{d_j(n)}{\sqrt{n}} \omega_j\left(n\left(\frac{\pi}{8d}\right)^{\frac{j}{2}}\right).$$

The relevance of these definitions is made clear in Lemma 2.2 where we show that for square-free odd integers $d$, and all integers $j \geq 1$,

$$L(\tfrac{1}{2}, \chi_{8d})^j = 2A_j(d).$$

From these approximate functional equations, we see that in order to evaluate the mollified moments we need asymptotic formulae for $\mathcal{S}_M(M(d)^j A_j(d); \Phi)$ (for $j = 1$, or 2). Further, we need good estimates for the remainder terms $\mathcal{S}_R(|M(d)^j A_j(d)|; \Phi)$ (for $j = 1$, or 2). In Section 3, we tackle the remainder terms and show that for "reasonable" mollifiers, their contribution is negligible.

PROPOSITION 1.1. *Suppose that $M(d)$ is as in* (1.3), *and that $\lambda(l) \ll l^{-1+\varepsilon}$. Then, for $j = 1, 2$,*

$$\mathcal{S}_R(|M(d)^j A_j(d)|; \Phi) \ll \frac{X^\varepsilon}{Y} + \frac{M^{\frac{j}{2}}}{X^{\frac{1}{2}-\varepsilon}}.$$

In Proposition 1.1 and throughout $\varepsilon$ denotes a small positive number. The reader should be warned that it might be a different $\varepsilon$ from line to line.

Next we evaluate $\mathcal{S}_M(M(d)A_1(d); \Phi)$. In fact, more generally we shall evaluate $\mathcal{S}_M(\left(\frac{8d}{l}\right)A_1(d); \Phi)$ where $l$ is any odd integer. Observe that

$$\mathcal{S}_M\left(\left(\frac{8d}{l}\right)A_1(d); \Phi\right) = \sum_{n=1}^{\infty} \frac{1}{\sqrt{n}} \mathcal{S}_M\left(\left(\frac{8d}{ln}\right); \Phi_n\right),$$

where $\Phi_n(t) = \Phi(t)\omega_1(n\sqrt{\pi}/\sqrt{8Xt})$. Now $\mathcal{S}_M(\left(\frac{8d}{ln}\right); \Phi_n)$ is essentially a character sum. Thus we may expect substantial cancellation here whenever $\left(\frac{\cdot}{ln}\right)$ is a nonprincipal character (i.e. $ln \neq \square$), and we may expect that the main term arises from the principal character terms $ln = \square$. Here, and throughout, we use the symbol $\square$ to denote square integers. In Section 4, we use the Pólya-Vinogradov inequality to make these heuristics precise, and establish Proposition 1.2.

PROPOSITION 1.2. *Write $l = l_1 l_2^2$ where $l_1$ and $l_2$ are odd with $l_1$ square-free. Then*

$$2\mathcal{S}_M\left(\left(\frac{8d}{l}\right)A_1(d);\Phi\right) = \frac{\hat{\Phi}(0)}{\zeta(2)\sqrt{l_1}}\frac{C}{g(l)}\left(\log\frac{\sqrt{X}}{l_1} + C_2 + \sum_{p|l}\frac{C_2(p)}{p}\log p\right)$$
$$+ O\left(\Phi_{(1)}\frac{l_1^{\frac{1}{2}+\varepsilon}Y}{X^{\frac{1}{2}-\varepsilon}} + \frac{\log X}{Y\sqrt{l_1}}\right),$$

*where*

$$C = \frac{1}{3}\prod_{p\geq 3}\left(1 - \frac{1}{p(p+1)}\right), \quad \text{and} \quad g(l) = \prod_{p|l}\left(\frac{p+1}{p}\right)\left(1 - \frac{1}{p(p+1)}\right).$$

*Lastly, $C_2$ is a constant depending only on $\Phi$ (it may be written as $C_3 + C_4\check{\Phi}'(0)/\check{\Phi}(0)$ for absolute constants $C_3$ and $C_4$) and $C_2(p) \ll 1$ for all $p$.*

Finally, it remains to handle $\mathcal{S}_M(|M(d)|^2 A_2(d);\Phi)$. Again, we treat the more general $\mathcal{S}_M((\frac{8d}{l})A_2(d);\Phi)$ where $l$ is any odd integer. As before, we may write

$$\mathcal{S}_M\left(\left(\frac{8d}{l}\right)A_2(d);\Phi\right) = \sum_{n=1}^{\infty}\frac{d(n)}{\sqrt{n}}\mathcal{S}_M\left(\left(\frac{8d}{ln}\right);F_n\right)$$

where $F_n(t) = \Phi(t)\omega_2(n\pi/8Xt)$. Again we expect that there is substantial cancellation in the character sum $\mathcal{S}_M((\frac{8d}{ln});F_n)$ when $ln \neq \square$, and that the main contribution comes from the $ln = \square$ terms. However, the simple Pólya-Vinogradov type argument of Section 4 is not enough to justify this; and, in fact, our expectation is wrong. There is an additional "off-diagonal" contribution to $\mathcal{S}_M((\frac{8d}{ln});F_n)$.

In Section 5, we develop a more delicate argument using Poisson summation to handle this (see Lemma 2.6 below). Roughly speaking, Poisson summation converts $\mathcal{S}_M((\frac{8d}{ln});F_n)$ into a sum of the form

$$\sum_k \left(\frac{k}{ln}\right)\tilde{F}_n\left(\frac{kX}{ln}\right)$$

where $\tilde{F}_n$ is essentially the Fourier transform of $F_n$. Now $\left(\frac{0}{ln}\right) = 1$ or $0$ depending on whether $ln$ is a square or not. So this term isolates the expected diagonal contribution of the terms $ln = \square$. The terms $k \neq 0$, or a $\square$ contribute a negligible amount because here $\left(\frac{k}{\cdot}\right)$ is a nonprincipal character. However, there is an additional contribution from the $k = \square$ terms which cannot be ignored. Evaluating this nondiagonal contribution forms the most subtle part of our argument, and we achieve this in Section 5.3. We note that these nondiagonal terms do not arise in the case of twisting a nonquadratic $L$-function $L(s,\psi)$ (as in [19]), because $\psi(\cdot)\left(\frac{k}{\cdot}\right)$ is nonprincipal for *all* $k \neq 0$.

For all integers $j \geq 0$ we define $\Lambda_j(n)$ to be the coefficient of $n^{-s}$ in the Dirichlet series expansion of $(-1)^j \zeta^{(j)}(s)/\zeta(s)$. Thus $\Lambda_0(1) = 1$, and $\Lambda_0(n) = 0$ for all $n \geq 2$; $\Lambda_1(n)$ is the usual von Mangoldt function $\Lambda(n)$. In general $\Lambda_j(n)$ is supported on integers having at most $j$ distinct prime factors, and $\Lambda_j(n) \ll_j (\log n)^j$.

PROPOSITION 1.3. *Write $l = l_1 l_2^2$ where $l_1$ and $l_2$ are odd and $l_1$ is square-free. Then*

$$2\mathcal{S}_M\left(\left(\frac{8d}{l}\right)A_2(d); \Phi\right)$$
$$= \frac{D\hat{\Phi}(0)}{36\zeta(2)} \frac{d(l_1)}{\sqrt{l_1}} \frac{l_1}{\sigma(l_1)h(l)} \left(\log^3\left(\frac{X}{l_1}\right) - 3\sum_{p|l_1} \log^2 p \log\left(\frac{X}{l_1}\right) + \mathcal{O}(l)\right)$$
$$+ \mathcal{R}(l) + O\left(\frac{l^\varepsilon X^\varepsilon}{\sqrt{l_1}Y} + \frac{l^\varepsilon X^\varepsilon}{(l_1 X)^{\frac{1}{4}}}\right),$$

*where $h$ is the multiplicative function defined on prime powers by*

$$h(p^k) = 1 + \frac{1}{p} + \frac{1}{p^2} - \frac{4}{p(p+1)}, \quad (k \geq 1), \quad D = \frac{1}{8}\prod_{p\geq 3}\left(1 - \frac{1}{p}\right)h(p),$$

*and*

$$\mathcal{O}(l) = \sum_{j,k=0}^{3} \sum_{m|l} \sum_{n|l_1} \frac{\Lambda_j(m)}{m} \frac{\Lambda_k(n)}{n} D(m,n) Q_{j,k}\left(\log \frac{X}{l_1}\right)$$
$$- 3\left(A + B\frac{\check{\Phi}'(0)}{\check{\Phi}(0)}\right) \sum_{p|l_1} \log^2 p$$

*where the $Q_{j,k}$ are polynomials of degree $\leq 2$ whose coefficients involve absolute constants and linear combinations of $\check{\Phi}^{(j)}(0)/\check{\Phi}(0)$ for $j = 1, 2, 3$; $A$ and $B$ are absolute constants; and $D(m,n) \ll 1$ uniformly for all $m$ and $n$. Lastly, $\mathcal{R}(l)$ is a remainder term bounded for each individual $l$ by*

$$|\mathcal{R}(l)| \ll \Phi_{(2)} \Phi_{(3)}^\varepsilon \frac{l^{\frac{1}{2}+\varepsilon} Y^{1+\varepsilon}}{X^{\frac{1}{2}-\varepsilon}},$$

*and bounded on average by*

$$\sum_{l=L}^{2L-1} |\mathcal{R}(l)| \ll \Phi_{(2)} \Phi_{(3)}^\varepsilon \frac{L^{1+\varepsilon} Y^{1+\varepsilon}}{X^{\frac{1}{2}-\varepsilon}}.$$

In Section 6 we choose our mollifier $M(d)$, and use Propositions 1.2 and 1.3 to complete the proof of Theorem 1. Our analysis there shows that an optimal mollifier of length $(\sqrt{X})^\theta$ leads to a proportion $\geq 1 - (\theta+1)^{-3} + o(1)$ for the nonvanishing of $L(\frac{1}{2}, \chi_{8d})$. Since Propositions 1.2 and 1.3 allow us to take a mollifier of length $X^{\frac{1}{2}-\varepsilon}$ we get the proportion $\frac{7}{8}$ of Theorem 1. If we believe that such results hold for mollifiers of arbitrary length (i.e. let $\theta \to \infty$), then we would get that $L(\frac{1}{2}, \chi_{8d}) \neq 0$ for almost all fundamental discriminants $8d$. We remark that Kowalski and Michel show that a mollifier of length $(\sqrt{q})^\theta$ (in their context of the rank of $J_0(q)$ [14]) leads to the same proportion $1 - (\theta + 1)^{-3}$ for the nonvanishing of $L'(f, \frac{1}{2})$. Curiously, this proportion also appears in a conditional result of J. B. Conrey, A. Ghosh, and S. M. Gonek [4] on simple zeros of $\zeta(s)$. They showed (assuming GRH) that a mollifier of length $T^\theta$ leads to a proportion $1 - (\theta + 1)^{-3}$ for the number of simple zeros of $\zeta(s)$ below height $T$. We gave earlier an explanation for the similarity between the Kowalski-Michel result and ours; it is unclear whether the similarity with this result of Conrey *et al.* is just a coincidence, or not.

We also note that using Proposition 1.3 with $l = 1$ we may deduce the following stronger form of Jutila's asymptotic formula (1.2).

COROLLARY 1.4. *There is a polynomial $Q$ of degree 3 such that*
$$\sum_{0 \leq d \leq X} L(\tfrac{1}{2}, \chi_{8d})^2 = XQ(\log X) + O(X^{\frac{5}{6}+\varepsilon}),$$
*where the sum is over fundamental discriminants $8d$.*

Corollary 1.4 should be compared with Heath-Brown's result on the fourth moment of $\zeta(s)$; see [8]. No doubt the remainder term in Corollary 1.4 can be refined; but we have not worried about optimizing it. Also one can calculate explicitly the coefficients of $Q(x)$ from our proof (compare Conrey [3]). Professor Heath-Brown has informed us that C. R. Guo (preprint) has obtained a result like Corollary 1.4 with a remainder term $O(X^{1-\frac{1}{1500}+\varepsilon})$.

While we cannot obtain an asymptotic formula for the fourth moment of $L(\frac{1}{2}, \chi_{8d})$, our methods enable us to evaluate the third moment.

THEOREM 2. *There is a polynomial $R$ of degree 6 such that*
$$\sum_{0 \leq d \leq X} L(\tfrac{1}{2}, \chi_{8d})^3 = XR(\log X) + O(X^{\frac{11}{12}+\varepsilon})$$
*where the sum is over fundamental discriminants $8d$.*

We shall merely sketch the proof of Theorem 2 in Section 7, since the details are very similar to the analysis carried out in other parts of this paper.

I am very grateful to Peter Sarnak for his constant encouragement and many helpful conversations. I also thank Brian Conrey and David Farmer for some useful conversations on the nature of the off-diagonal contribution to Proposition 1.3. Lastly I am grateful to the referee for a careful reading of the manuscript and some valuable suggestions.

## 2. Preliminaries

2.1. *Approximate functional equations.* We first prove some properties of the functions $\omega_j(\xi)$ defined in (1.5).

LEMMA 2.1. *The functions $\omega_j(\xi)$ are real-valued and smooth on $[0, \infty)$. For $\xi$ near $0$ they satisfy*

$$\omega_j(\xi) = 1 + O(\xi^{\frac{1}{2}-\varepsilon}),$$

*and for large $\xi$ and any integer $\nu$,*

$$\omega_j^{(\nu)}(\xi) \ll_{\nu,j} \xi^{2\nu+2} \exp(-j\xi^{\frac{2}{j}}) \ll_{\nu,j} \exp(-\tfrac{j}{2}\xi^{\frac{2}{j}}).$$

*Proof.* By pairing together the $s$ and $\bar{s}$ values of the integrand in (1.5), we see that $\omega_j(\xi)$ is real-valued. Further the $\nu$-th derivative of $\omega_j(\xi)$ is plainly

$$(2.1) \qquad \frac{(-1)^\nu}{2\pi i} \int_{(c)} \left(\frac{\Gamma(\frac{s}{2}+\frac{1}{4})}{\Gamma(\frac{1}{4})}\right)^j s(s+1)\ldots(s+\nu-1)\xi^{-s}\frac{ds}{s}.$$

Thus $\omega_j(\xi)$ is smooth.

Move the line of integration in (1.5) to the line from $-\frac{1}{2}+\varepsilon - i\infty$ to $-\frac{1}{2}+\varepsilon+i\infty$. The pole at $s=0$ leaves the residue $1$, and the integral on the new line is plainly $\ll \xi^{\frac{1}{2}-\varepsilon}$. Thus $\omega_j(\xi) = 1 + O(\xi^{\frac{1}{2}-\varepsilon})$, as desired.

To prove the last estimate of the lemma, we may suppose that $\xi^{\frac{2}{j}} \geq \nu+2$. Since $|\Gamma(x+iy)| \leq \Gamma(x)$ for $x \geq 1$, and $s\Gamma(s) = \Gamma(s+1)$ we obtain that (2.1) is (here $c>0$ is arbitrary)

$$\ll_\nu \Gamma(\tfrac{c}{2}+\tfrac{1}{4}+\nu+1)^j \xi^{-c} \int_{(c)} \frac{|ds|}{|s|^2} \ll_\nu \Gamma(\tfrac{c}{2}+\nu+\tfrac{5}{4})^j \frac{\xi^{-c}}{c}.$$

By Stirling's formula this is

$$\ll_\nu \left(\frac{c+2\nu+2}{2e}\right)^{\frac{j}{2}(c+2\nu+2)} \frac{\xi^{-c}}{c}.$$

With $c = 2\xi^{\frac{2}{j}} - 2\nu - 2(\geq 2)$ above, the desired estimate follows.

Recall from (1.6) the definition of $A_j(d)$.

LEMMA 2.2. *Suppose that $d$ is an odd, positive, square-free number. Then, for all integers $j \geq 1$,*

$$L(\tfrac{1}{2}, \chi_{8d})^j = 2A_j(d).$$

*Proof.* For some $c > \frac{1}{2}$ consider

$$(2.2) \qquad \frac{1}{2\pi i} \int_{(c)} \left(\frac{\Gamma(\frac{s}{2}+\frac{1}{4})}{\Gamma(\frac{1}{4})}\right)^j L(s+\tfrac{1}{2}, \chi_{8d})^j \left(\frac{8d}{\pi}\right)^{j\frac{s}{2}} \frac{ds}{s}.$$

Expanding $L(s+\frac{1}{2}, \chi_{8d})^j$ into its Dirichlet series we see that this equals $A_j(d)$. We now evaluate (2.2) differently by moving the line of integration to the Re $(s) = -\frac{1}{8}$ line. The pole at $s=0$ leaves the residue $L(\frac{1}{2}, \chi_{8d})^j$. Thus (2.2) equals

$$(2.3) \qquad L(\tfrac{1}{2}, \chi_{8d})^j + \frac{1}{2\pi i}\int_{(-\frac{1}{8})} \left(\frac{\Gamma(\frac{s}{2}+\frac{1}{4})}{\Gamma(\frac{1}{4})}\right)^j L(\tfrac{1}{2}+s, \chi_{8d})^j \left(\frac{8d}{\pi}\right)^{j\frac{s}{2}} \frac{ds}{s}.$$

Recall from [6, Chap. 9] the functional equation for $L(\frac{1}{2}+s, \chi_{8d})$:

$$\left(\frac{8d}{\pi}\right)^{\frac{s}{2}} \Gamma(\tfrac{s}{2}+\tfrac{1}{4}) L(\tfrac{1}{2}+s, \chi_{8d}) = \frac{\tau(\chi_{8d})}{\sqrt{8d}} \left(\frac{8d}{\pi}\right)^{-\frac{s}{2}} \Gamma(\tfrac{1}{4}-\tfrac{s}{2}) L(\tfrac{1}{2}-s, \chi_{8d}).$$

Here $\tau(\chi_{8d})$ is the Gauss sum of $\chi_{8d}$ (mod $8d$). Since $8d$ is a fundamental discriminant we note that $\tau(\chi_{8d}) = \sqrt{8d}$ (see [6, Chap. 2]). From this it follows that the integral in (2.3) equals

$$\frac{1}{2\pi i}\int_{(-\frac{1}{8})} \left(\frac{\Gamma(\frac{1}{4}-\frac{s}{2})}{\Gamma(\frac{1}{4})}\right)^j L(\tfrac{1}{2}-s, \chi_{8d})^j \left(\frac{8d}{\pi}\right)^{-j\frac{s}{2}} \frac{ds}{s}.$$

Replacing $s$ by $-s$ we see that the above equals $-A_j(d)$; and this gives the lemma.

2.2. *On Gauss-type sums.* Let $n$ be an odd integer. We define for all integers $k$

$$G_k(n) = \left(\frac{1-i}{2} + \left(\frac{-1}{n}\right)\frac{1+i}{2}\right) \sum_{a \pmod n} \left(\frac{a}{n}\right) e\left(\frac{ak}{n}\right),$$

and put

$$\tau_k(n) = \sum_{a \pmod n} \left(\frac{a}{n}\right) e\left(\frac{ak}{n}\right) = \left(\frac{1+i}{2} + \left(\frac{-1}{n}\right)\frac{1-i}{2}\right) G_k(n).$$

If $n$ is square-free then $\left(\frac{\cdot}{n}\right)$ is a primitive character with conductor $n$. Here it is easy to see that $G_k(n) = \left(\frac{k}{n}\right)\sqrt{n}$. For our later work, we require knowledge of $G_k(n)$ for all odd $n$. In the next lemma we show how this may be attained.

LEMMA 2.3. (i) (*Multiplicativity*) *Suppose $m$ and $n$ are coprime odd integers. Then $G_k(mn) = G_k(m)G_k(n)$.*
(ii) *Suppose $p^\alpha$ is the largest power of $p$ dividing $k$. (If $k=0$ then set $\alpha = \infty$.) Then for $\beta \geq 1$*

$$G_k(p^\beta) = \begin{cases} 0 & \text{if } \beta \leq \alpha \text{ is odd}, \\ \varphi(p^\beta) & \text{if } \beta \leq \alpha \text{ is even}, \\ -p^\alpha & \text{if } \beta = \alpha+1 \text{ is even}, \\ \left(\frac{kp^{-\alpha}}{p}\right) p^\alpha \sqrt{p} & \text{if } \beta = \alpha+1 \text{ is odd} \\ 0 & \text{if } \beta \geq \alpha+2. \end{cases}$$

*Proof.* Using the Chinese remainder theorem we may write $a \pmod{mn}$ as $bn + cm$ with $b \pmod m$ and $c \pmod n$. This shows that $\tau_k(mn) = \left(\frac{m}{n}\right)\left(\frac{n}{m}\right)\tau_k(m)\tau_k(n)$. To show (i) we need only check that

$$\frac{1-i}{2} + \left(\frac{-1}{mn}\right)\frac{1+i}{2} = \left(\frac{1-i}{2} + \left(\frac{-1}{m}\right)\frac{1+i}{2}\right)\left(\frac{1-i}{2} + \left(\frac{-1}{n}\right)\frac{1+i}{2}\right)\left(\frac{m}{n}\right)\left(\frac{n}{m}\right);$$

this holds by quadratic reciprocity.

If $\beta = \alpha + 1$ then

$$\sum_{a \pmod{p^\beta}} \left(\frac{a}{p^\beta}\right) e\left(\frac{ak}{p^\beta}\right) = \sum_{l \pmod p} \left(\frac{l}{p^\beta}\right) \sum_{b \pmod{p^{\beta-1}}} e\left(\frac{(bp+l)k}{p^\beta}\right)$$

$$= p^{\beta-1} \sum_{l \pmod p} \left(\frac{l}{p^\beta}\right) e\left(\frac{lk}{p^\beta}\right).$$

If $\beta$ is even then the last sum above is $-1$ and if $\beta$ is odd the last sum above is, from knowledge of the usual Gauss sum (see [6, Chap. 2]),

$$\sum_{l \pmod p} \left(\frac{l}{p}\right) e\left(\frac{l(kp^{-\alpha})}{p}\right) = \left(\frac{kp^{-\alpha}}{p}\right) \times \begin{cases} \sqrt{p} & \text{if } p \equiv 1 \pmod 4 \\ i\sqrt{p} & \text{if } p \equiv 3 \pmod 4. \end{cases}$$

These calculations show the third and fourth cases of (ii). The other cases are left as easy exercises for the reader.

2.3. *Lemmas for estimating character sums.* We collect here two lemmas that will be very useful in bounding the character sums that arise below. These are consequences of a recent large sieve result for real characters due to D. R. Heath-Brown [9].

LEMMA 2.4. *Let $N$ and $Q$ be positive integers and let $a_1, \ldots, a_N$ be arbitrary complex numbers. Let $S(Q)$ denote the set of real, primitive characters $\chi$ with conductor $\leq Q$. Then*

$$\sum_{\chi \in S(Q)} \left|\sum_{n \leq N} a_n \chi(n)\right|^2 \ll_\varepsilon (QN)^\varepsilon (Q+N) \sum_{n_1 n_2 = \square} |a_{n_1} a_{n_2}|,$$

*for any $\varepsilon > 0$. Let $M$ be any positive integer, and for each $|m| \leq M$ write $4m = m_1 m_2^2$ where $m_1$ is a fundamental discriminant, and $m_2$ is positive. Suppose the sequence $a_n$ satisfies $|a_n| \ll n^\varepsilon$. Then*

$$\sum_{|m| \leq M} \frac{1}{m_2} \left|\sum_{n \leq N} a_n \left(\frac{m}{n}\right)\right|^2 \ll (MN)^\varepsilon N(M+N).$$

*Proof.* The first assertion is Corollary 2 of Heath-Brown [9]. Using this result, we see that the second quantity to be bounded is

$$\ll \sum_{m_2 \leq 2\sqrt{M}} \frac{1}{m_2} \sum_{\chi \in S(M/m_2^2)} \left| \sum_{n \leq N} a_n\left(\frac{m_2^2}{n}\right) \chi(n) \right|^2$$

$$\ll \sum_{m_2 \leq 2\sqrt{M}} \frac{1}{m_2} (NM)^\varepsilon N\left(N + \frac{M}{m_2^2}\right),$$

and the result follows.

LEMMA 2.5. *Let $S(Q)$ be as in Lemma 2.4, and suppose $\sigma + it$ is a complex number with $\sigma \geq \frac{1}{2}$. Then*

$$\sum_{\chi \in S(Q)} |L(\sigma + it, \chi)|^4 \ll Q^{1+\varepsilon}(1+|t|)^{1+\varepsilon},$$

*and*

$$\sum_{\chi \in S(Q)} |L(\sigma + it, \chi)|^2 \ll Q^{1+\varepsilon}(1+|t|)^{\frac{1}{2}+\varepsilon}.$$

*Proof.* The fourth moment estimate is in Theorem 2 of Heath-Brown [9]. The second moment estimate follows from this by Cauchy's inequality.

2.4 *Poisson summation.* For a Schwarz class function $F$ we define

$$\tilde{F}(\xi) = \frac{1+i}{2}\hat{F}(\xi) + \frac{1-i}{2}\hat{F}(-\xi) = \int_{-\infty}^{\infty} (\cos(2\pi\xi x) + \sin(2\pi\xi x)) F(x) dx.$$

LEMMA 2.6. *Let $F$ be a nonnegative, smooth function supported in $(1,2)$. For any odd integer $n$,*

$$\mathcal{S}_M\left(\left(\frac{d}{n}\right); F\right) = \frac{1}{2n}\left(\frac{2}{n}\right) \sum_{\substack{\alpha \leq Y \\ (\alpha, 2n)=1}} \frac{\mu(\alpha)}{\alpha^2} \sum_k (-1)^k G_k(n) \tilde{F}\left(\frac{kX}{2\alpha^2 n}\right).$$

*Proof.* First note that

$$(2.4) \quad \sum_{\substack{d \\ (d,2)=1}} M_Y(d)\left(\frac{d}{n}\right) F\left(\frac{d}{X}\right) = \sum_{\substack{\alpha \leq Y \\ (\alpha, 2n)=1}} \mu(\alpha) \sum_{\substack{d \\ (d,2)=1}} \left(\frac{d}{n}\right) F\left(\frac{d\alpha^2}{X}\right).$$

Next observe that

$$(2.5) \quad \sum_{d \text{ odd}} \left(\frac{d}{n}\right) F\left(\frac{d\alpha^2}{X}\right) = \sum_d \left(\frac{d}{n}\right) F\left(\frac{d\alpha^2}{X}\right) - \left(\frac{2}{n}\right) \sum_d \left(\frac{d}{n}\right) F\left(\frac{2d\alpha^2}{X}\right).$$

Splitting the sum over $d$ below according to the residue classes $\pmod{n}$ and using the Poisson summation formula we derive (for $a = 1$, or 2)

$$\sum_d \left(\frac{d}{n}\right) F\left(\frac{ad\alpha^2}{X}\right) = \sum_{b \pmod{n}} \left(\frac{b}{n}\right) \sum_d F\left(\frac{a\alpha^2(nd+b)}{X}\right)$$

$$= \frac{X}{na\alpha^2} \sum_{b \pmod{n}} \left(\frac{b}{n}\right) \sum_k \hat{F}\left(\frac{kX}{an\alpha^2}\right) e\left(\frac{kb}{n}\right)$$

$$= \frac{X}{na\alpha^2} \sum_k \hat{F}\left(\frac{kX}{a\alpha^2 n}\right) \tau_k(n).$$

Writing $\tau_k$ in terms of $G_k$, using the relation $G_k(n) = \left(\frac{-1}{n}\right) G_{-k}(n)$, and recombining the $k$ and $-k$ terms, we obtain that the above is

$$\frac{X}{na\alpha^2} \sum_k G_k(n) \tilde{F}\left(\frac{kX}{a\alpha^2 n}\right).$$

Substituting this in the right-hand side of (2.5) we see that (using $G_k(n) = \left(\frac{2}{n}\right) G_{2k}(n)$)

$$(2.5) = \frac{X}{2n\alpha^2} \left(\frac{2}{n}\right) \sum_k (-1)^k G_k(n) \tilde{F}\left(\frac{kX}{2n\alpha^2}\right).$$

Substituting this in (2.4) we get the lemma.

## 3. Proof of Proposition 1.1

Observe that $R_Y(d) = 0$ unless $d = l^2 m$ where $m$ is square-free and $l > Y$. Further, note that $|R_Y(d)| \leq \sum_{k|d} 1 \ll d^\varepsilon$. Hence

$$\mathcal{S}_R(|M(d)^j A_j(d)|; \Phi) \ll X^{-1+\varepsilon} \sum_{\substack{Y < l \leq \sqrt{2X} \\ (l,2)=1}} \sum_{X/l^2 \leq m \leq 2X/l^2}^{\flat} |M(l^2 m)^j A_j(l^2 m)|,$$

where the $\flat$ on the sum over $m$ indicates that $m$ is odd and square-free, and $j = 1$, or 2. By Cauchy's inequality the above is
(3.1)

$$\ll X^{-1+\varepsilon} \sum_{\substack{Y < l \leq \sqrt{2X} \\ (l,2)=1}} \left( \sum_{X/l^2 \leq m \leq 2X/l^2}^{\flat} |M(l^2 m)|^{2j} \right)^{\frac{1}{2}} \left( \sum_{X/l^2 \leq m \leq 2X/l^2}^{\flat} |A_j(l^2 m)|^2 \right)^{\frac{1}{2}}.$$

Write $\lambda_1(n) = \lambda(n)$ and $\lambda_2(n) = \sum_{ab=n, a,b \leq M} \lambda(a)\lambda(b)$. Note that $|\lambda_j(n)| \ll n^{-1+\varepsilon}$ and that $M(d)^j = \sum_{n \leq M^j} \lambda_j(n) \sqrt{n} \left(\frac{8d}{n}\right)$. Hence, by Lemma 2.4,

$$
\text{(3.2)} \quad \sum_{X/l^2 \leq m \leq 2X/l^2}^{\flat} |M(l^2 m)|^{2j} = \sum_{X/l^2 \leq m \leq 2X/l^2}^{\flat} \left| \sum_{n \leq M^j} \lambda_j(n) \sqrt{n} \left(\frac{l^2}{n}\right) \left(\frac{8m}{n}\right) \right|^2
$$

$$
\ll X^{\varepsilon} \left(\frac{X}{l^2} + M^j\right) \sum_{\substack{n_1, n_2 \leq M^j \\ n_1 n_2 = \square}} |\lambda_j(n_1) \lambda_j(n_2)| \sqrt{n_1 n_2}
$$

$$
\ll X^{\varepsilon} \left(\frac{X}{l^2} + M^j\right) \sum_{\substack{n_1, n_2 \leq M^j \\ n_1 n_2 = \square}} \frac{1}{\sqrt{n_1 n_2}}
$$

$$
\ll X^{\varepsilon} \left(\frac{X}{l^2} + M^j\right) \sum_{a \leq M^{2j}} \frac{d(a^2)}{a} \ll X^{\varepsilon} \left(\frac{X}{l^2} + M^j\right).
$$

Now observe that for any $c > \frac{1}{2}$,

$$
\text{(3.3)} \quad A_j(l^2 m) = \sum_{n=1}^{\infty} \frac{d_j(n)}{\sqrt{n}} \left(\frac{8l^2 m}{n}\right) \omega_j \left( n \left(\frac{\pi}{8l^2 m}\right)^{\frac{j}{2}} \right)
$$

$$
= \frac{1}{2\pi i} \int_{(c)} \left( \frac{\Gamma(\frac{s}{2} + \frac{1}{4})}{\Gamma(\frac{1}{4})} \right)^j \left( \frac{8l^2 m}{\pi} \right)^{s\frac{j}{2}} \sum_{n=1}^{\infty} \frac{d_j(n)}{n^{s+\frac{1}{2}}} \left(\frac{8l^2 m}{n}\right) \frac{ds}{s}.
$$

Plainly

$$
\text{(3.4)} \quad \sum_{n=1}^{\infty} \frac{d_j(n)}{n^{s+\frac{1}{2}}} \left(\frac{8l^2 m}{n}\right) = L(\tfrac{1}{2} + s, \chi_{8m})^j \mathcal{E}(s, l)^j
$$

where

$$
\mathcal{E}(s, l) = \prod_{p | l} \left( 1 - \frac{1}{p^{s+\frac{1}{2}}} \left(\frac{8m}{p}\right) \right).
$$

Since $\chi_{8m}$ is nonprincipal, it follows that the left side of (3.4) is analytic for all $s$.

Hence we may move the line of integration in (3.3) to the line from $1/\log X - i\infty$ to $1/\log X + i\infty$. This gives

$$
|A_j(l^2 m)| \ll \int_{(\frac{1}{\log X})} |\Gamma(\tfrac{s}{2} + \tfrac{1}{4})|^j |L(\tfrac{1}{2} + s, \chi_{8m})|^j |\mathcal{E}(s, l)|^j \frac{|ds|}{|s|}.
$$

Plainly $|\mathcal{E}(s, l)| \leq \prod_{p|l}(1 + 1/\sqrt{p}) \ll l^{\varepsilon} \ll X^{\varepsilon}$, and note that

$$
\int_{(\frac{1}{\log X})} |\Gamma(\tfrac{s}{2} + \tfrac{1}{4})|^j \frac{|ds|}{|s|^2} \ll X^{\varepsilon}.
$$

Using these estimates and Cauchy's inequality, we deduce

$$|A_j(l^2m)|^2 \ll X^\varepsilon \int_{(\frac{1}{\log X})} |\Gamma(\tfrac{s}{2}+\tfrac{1}{4})|^j |L(\tfrac{1}{2}+s,\chi_{8m})|^{2j} |ds|.$$

Summing this over $m$ and using Lemma 2.5, we obtain
(3.5)
$$\sum_{X/l^2 \leq m \leq 2X/l^2}^{\flat} |A_j(l^2m)|^2 \ll \frac{X^{1+\varepsilon}}{l^2} \int_{(\frac{1}{\log X})} |\Gamma(\tfrac{s}{2}+\tfrac{1}{4})|^2 (1+|s|)^{1+\varepsilon} |ds| \ll \frac{X^{1+\varepsilon}}{l^2}.$$

Proposition 1.1 follows upon combination of (3.1) with (3.2) and (3.5).

## 4. Proof of Proposition 1.2

Observe that

(4.1) $$\mathcal{S}_M\left(\left(\frac{8d}{l}\right)A_1(d);\Phi\right) = \sum_{n=1}^\infty \frac{1}{\sqrt{n}} \mathcal{S}_M\left(\left(\frac{8d}{ln}\right);\Phi_n\right),$$

where

$$\Phi_n(t) = \Phi(t)\omega_1\left(\frac{n\sqrt{\pi}}{\sqrt{8Xt}}\right).$$

LEMMA 4.1. *If $ln \neq \square$ then*

$$\mathcal{S}_M\left(\left(\frac{8d}{ln}\right);\Phi_n\right) \ll \Phi_{(1)} \frac{Y}{X} \sqrt{ln}\log(ln)\exp\left(-\frac{n}{10X^{\frac{1}{2}}}\right).$$

*If $ln = \square$ then*

$$\mathcal{S}_M\left(\left(\frac{8d}{ln}\right);\Phi_n\right) = \left(\frac{8}{ln}\right)\frac{\hat{\Phi}_n(0)}{\zeta(2)} \prod_{p|2ln}\left(\frac{p}{p+1}\right)\left(1+O\left(\frac{1}{Y}\right)\right)$$
$$+ O\left(\Phi_{(1)}\frac{Yl^\varepsilon n^\varepsilon}{X}\exp\left(-\frac{n}{10X^{\frac{1}{2}}}\right)\right).$$

*Proof.* Note that $\left(\frac{d}{4ln}\right) = \left(\frac{d}{ln}\right)$ if $d$ is odd and is 0 otherwise. Thus we seek to bound (or evaluate)

(4.2) $$\mathcal{S}_M\left(\left(\frac{8d}{ln}\right);\Phi_n\right) = \frac{1}{X}\sum_{\substack{\alpha \leq Y \\ \alpha \text{ odd}}} \mu(\alpha)\left(\frac{8\alpha^2}{ln}\right)\sum_d \left(\frac{d}{4ln}\right)\Phi_n\left(\frac{d\alpha^2}{X}\right).$$

If $ln \neq \square$, $\left(\frac{\cdot}{4ln}\right)$ is a nonprincipal character to the modulus $4ln$. Hence by the Pólya-Vinogradov inequality

$$\sum_{d \leq x}\left(\frac{d}{4ln}\right) \ll \sqrt{ln}\log(4ln)$$

for all $x$. By partial summation it follows that
$$\sum_d \left(\frac{d}{4ln}\right) \Phi_n\left(\frac{d\alpha^2}{X}\right) \ll \sqrt{ln} \log(4ln) \int_0^\infty |\Phi_n'(t)|dt.$$

Now, by Lemma 2.1,

(4.3)
$$\int_0^\infty |\Phi_n'(t)|dt \leq \int_1^2 \left(|\Phi'(t)|\omega_1\left(\frac{n\sqrt{\pi}}{\sqrt{8Xt}}\right) + \Phi(t)\left|\omega_1'\left(\frac{n\sqrt{\pi}}{\sqrt{8Xt}}\right)\right|\frac{n\sqrt{\pi}}{2\sqrt{8Xt^3}}\right) dt$$
$$\ll \Phi_{(1)} \exp\left(-\frac{n}{10X^{\frac{1}{2}}}\right).$$

Using these estimates in (4.2), we obtain the first bound of the lemma.

If $ln = \square$ then $\left(\frac{d}{4ln}\right) = 1$ if $d$ is coprime to $2ln$, and $0$ otherwise. Hence
$$\sum_{d \leq x} \left(\frac{d}{4ln}\right) = \frac{\varphi(2ln)}{2ln}x + O((ln)^\varepsilon),$$

and so by partial summation and (4.3) we get
$$\sum_d \left(\frac{d}{4ln}\right)\Phi_n\left(\frac{d\alpha^2}{X}\right) = \frac{\varphi(2ln)}{2ln}\frac{X}{\alpha^2}\hat{\Phi}_n(0) + O\left(\Phi_{(1)}(ln)^\varepsilon \exp\left(-\frac{n}{10X^{\frac{1}{2}}}\right)\right).$$

We use this in (4.2) and observe that

(4.4) $$\sum_{\substack{\alpha \leq Y \\ (\alpha, 2ln)=1}} \frac{\mu(\alpha)}{\alpha^2} = \frac{1}{\zeta(2)} \prod_{p|2ln}\left(1 - \frac{1}{p^2}\right)^{-1}\left(1 + O\left(\frac{1}{Y}\right)\right).$$

This proves the second part of the lemma.

Using Lemma 4.1 in (4.1), we obtain

(4.5) $$\mathcal{S}_M\left(\left(\frac{8d}{l}\right)A_1(d); \Phi\right) = M(1 + O(Y^{-1})) + R,$$

where
$$M = \frac{1}{\zeta(2)} \sum_{\substack{n=1 \\ ln=\square}}^\infty \frac{1}{\sqrt{n}}\left(\frac{8}{ln}\right)\prod_{p|2ln}\left(\frac{p}{p+1}\right)\hat{\Phi}_n(0),$$

and

(4.6) $$R \ll \Phi_{(1)}\frac{Y}{X}\sum_{n=1}^\infty l^{\frac{1}{2}+\varepsilon}n^\varepsilon \exp\left(-\frac{n}{10\sqrt{X}}\right) \ll \Phi_{(1)}\frac{l^{\frac{1}{2}+\varepsilon}Y}{X^{\frac{1}{2}-\varepsilon}}.$$

We now focus on evaluating $M$. Recall that $l = l_1 l_2^2$ where $l_1$ and $l_2$ are odd and $l_1$ is square-free. Thus the condition $ln = \square$ is equivalent to $n = l_1 m^2$ for some integer $m$. Hence

$$(4.7) \quad M = \frac{1}{\sqrt{l_1}\zeta(2)} \sum_{\substack{m=1 \\ m \text{ odd}}}^{\infty} \frac{1}{m} \prod_{p|2lm} \left(\frac{p}{p+1}\right) \hat{\Phi}_{l_1 m^2}(0)$$

$$= \frac{1}{\sqrt{l_1}\zeta(2)} \int_1^2 \Phi(t) \sum_{\substack{m=1 \\ m \text{ odd}}}^{\infty} \frac{1}{m} \prod_{p|2lm} \left(\frac{p}{p+1}\right) \omega_1\left(\frac{l_1 m^2 \sqrt{\pi}}{\sqrt{8Xt}}\right) dt.$$

Write $\prod_{p|2lm}\left(\frac{p}{p+1}\right) = \sum_{d|2lm} \frac{\mu(d)}{\sigma(d)}$ where $\sigma(d)$ is the sum of the divisors of $d$. An elementary argument based on swapping the sums over $d$ and $m$ (which we leave to the reader) shows that

$$\sum_{\substack{m \le x \\ m \text{ odd}}} \frac{1}{m} \prod_{p|2lm} \left(\frac{p}{p+1}\right) = \frac{C}{g(l)}\left(\log x + C_0 + \sum_{p|l} \frac{C_0(p)}{p} \log p\right) + O\left(\frac{d(l)\log x}{x}\right),$$

where $C$ and $g(l)$ are as defined in the statement of Proposition 1.2, $C_0$ is an absolute constant, and $C_0(p) \ll 1$ for all primes $p$. Hence for any $1 \le t \le 2$ we have by partial summation (using Lemma 2.1)

$$\sum_{\substack{m=1 \\ m \text{ odd}}}^{\infty} \frac{1}{m} \prod_{p|2lm} \left(\frac{p}{p+1}\right) \omega_1\left(\frac{l_1 m^2 \sqrt{\pi}}{\sqrt{8Xt}}\right) = \frac{C}{g(l)}\left(\log \frac{X^{\frac{1}{4}}t^{\frac{1}{4}}}{l_1^{\frac{1}{2}}} + C_1 + \sum_{p|l} \frac{C_1(p)}{p} \log p\right)$$

$$+ O\left(\frac{d(l)l_1^{\frac{1}{2}}}{X^{\frac{1}{4}-\varepsilon}}\right),$$

where $C_1$ is an absolute constant and $C_1(p) \ll 1$ for all $p$. We use this expression in (4.7) to evaluate $M$. Combining this with (4.5) and (4.6), we see that Proposition 1.2 follows.

## 5. Proof of Proposition 1.3

Observe that

$$(5.1) \quad \mathcal{S}_M\left(\left(\frac{8d}{l}\right)A_2(d); \Phi\right) = \sum_{n=1}^{\infty} \frac{d(n)}{\sqrt{n}} \mathcal{S}_M\left(\left(\frac{8d}{ln}\right); F_n\right),$$

where

$$F_n(t) = \Phi(t)\omega_2\left(\frac{n\pi}{8Xt}\right).$$

Using Poisson summation, Lemma 2.6 above, we obtain
(5.2)
$$\mathcal{S}_M\left(\left(\frac{8d}{ln}\right); F_n\right) = \frac{1}{2ln}\left(\frac{16}{ln}\right) \sum_{\substack{\alpha \le Y \\ (\alpha,2ln)=1}} \frac{\mu(\alpha)}{\alpha^2} \sum_{k=-\infty}^{\infty} (-1)^k G_k(n) \tilde{F}_n\left(\frac{kX}{2\alpha^2 ln}\right).$$

Using this in (5.1) we deduce
$$\mathcal{S}_M\left(\left(\frac{8d}{l}\right)A_2(d);\Phi\right) = \mathcal{P}_1(l) + \mathcal{R}_0(l),$$
where $\mathcal{P}_1(l)$ is the main principal term (arising from the $k = 0$ term in (5.2)), and $\mathcal{R}_0(l)$ includes all the nonzero terms $k$ in (5.2). Thus
$$\mathcal{P}_1(l) = \frac{1}{2l}\sum_{n=1}^{\infty}\frac{d(n)}{n^{\frac{3}{2}}}\left(\frac{16}{ln}\right)\sum_{\substack{\alpha \leq Y \\ (\alpha, 2ln)=1}}\frac{\mu(\alpha)}{\alpha^2}G_0(ln)\tilde{F}_n(0),$$
and
$$(5.3)$$
$$\mathcal{R}_0(l) = \frac{1}{2l}\sum_{n=1}^{\infty}\frac{d(n)}{n^{\frac{3}{2}}}\left(\frac{16}{ln}\right)\sum_{\substack{\alpha \leq Y \\ (\alpha, 2ln)=1}}\frac{\mu(\alpha)}{\alpha^2}\sum_{\substack{k=-\infty \\ k \neq 0}}^{\infty}(-1)^k G_k(ln)\tilde{F}_n\left(\frac{kX}{2\alpha^2 ln}\right).$$

We shall compute the main principal contribution $\mathcal{P}_1(l)$ in Section 5.1 below. In Section 5.2, we separate $\mathcal{R}_0(l)$ into a second main term (essentially arising from the $k = \square$ terms in (5.2)), $\mathcal{P}_2(l)$, and a remainder term $\mathcal{R}(l)$. The evaluation of the $\mathcal{P}_2(l)$ contribution is quite subtle and forms the focus of our attention in Section 5.3. The remainder terms $\mathcal{R}(l)$ are relatively straightforward, and we bound their effect precisely in Section 5.4 below.

5.1. *The principal $\mathcal{P}_1(l)$ contribution.* Note that $\tilde{F}_n(0) = \hat{F}_n(0)$ and that $G_0(ln) = \varphi(ln)$ if $ln = \square$ and $G_0(ln) = 0$ otherwise. Using (4.4) and these observations we get
$$\mathcal{P}_1(l) = \frac{1 + O(Y^{-1})}{\zeta(2)}\sum_{\substack{n=1 \\ ln=\square}}^{\infty}\frac{d(n)}{n^{\frac{1}{2}}}\left(\frac{16}{ln}\right)\prod_{p|2ln}\left(\frac{p}{p+1}\right)\hat{F}_n(0).$$

Recall that $l = l_1 l_2^2$ where $l_1$ and $l_2$ are odd, and $l_1$ is square-free. The condition that $ln = \square$ is thus equivalent to $n = l_1 m^2$ for some integer $m$. Hence
$$\mathcal{P}_1(l) = \frac{1 + O(Y^{-1})}{\zeta(2)\sqrt{l_1}}\sum_{\substack{m=1 \\ m \text{ odd}}}^{\infty}\frac{d(l_1 m^2)}{m}\prod_{p|2lm}\left(\frac{p}{p+1}\right)\hat{F}_{l_1 m^2}(0).$$

For any $c > 0$,
$$\hat{F}_{l_1 m^2}(0) = \int_0^{\infty}\Phi(t)\omega_2\left(\frac{l_1 m^2 \pi}{8Xt}\right)dt$$
$$= \frac{1}{2\pi i}\int_{(c)}\left(\frac{\Gamma(\frac{s}{2} + \frac{1}{4})}{\Gamma(\frac{1}{4})}\right)^2\left(\frac{8X}{l_1 m^2 \pi}\right)^s\left(\int_0^{\infty}\Phi(t)t^s dt\right)\frac{ds}{s}$$
$$= \frac{1}{2\pi i}\int_{(c)}\left(\frac{\Gamma(\frac{s}{2} + \frac{1}{4})}{\Gamma(\frac{1}{4})}\right)^2\left(\frac{8X}{l_1 m^2 \pi}\right)^s\check{\Phi}(s)\frac{ds}{s}.$$

Thus
$$\mathcal{P}_1(l) = \frac{2}{3}\frac{1+O(Y^{-1})}{\zeta(2)\sqrt{l_1}}\frac{1}{2\pi i}\int_{(c)}\left(\frac{\Gamma(\frac{s}{2}+\frac{1}{4})}{\Gamma(\frac{1}{4})}\right)^2\left(\frac{8X}{l_1\pi}\right)^s\check{\Phi}(s)$$
$$\times \sum_{\substack{m=1 \\ m \text{ odd}}}^{\infty} \frac{d(l_1 m^2)}{m^{1+2s}} \prod_{p|lm}\left(\frac{p}{p+1}\right)\frac{ds}{s}.$$

LEMMA 5.1. *Suppose $l = l_1 l_2^2$ is as above. Then for $\mathrm{Re}(s) > 1$*
$$\sum_{\substack{m=1 \\ m \text{ odd}}}^{\infty} \frac{d(l_1 m^2)}{m^s} \prod_{p|lm}\left(\frac{p}{p+1}\right) = d(l_1)\zeta(s)^3 \eta(s;l)$$

where $\eta(s;l) = \prod_p \eta_p(s;l)$ with $\eta_2(s;l) = (1-2^{-s})^3$ and for $p \geq 3$,

$$\eta_p(s;l) = \begin{cases} 1 - \frac{4}{p^s(p+1)}\left(1-\frac{1}{p^s}\right) + \frac{1}{p^s(p+1)} - \frac{1}{p^{2s}} - \frac{1}{p^{3s}(p+1)} & \text{if } p \nmid l \\ \left(\frac{p}{p+1}\right)\left(1 - \frac{1}{p^s}\right) & \text{if } p|l_1 \\ \left(\frac{p}{p+1}\right)\left(1 - \frac{1}{p^{2s}}\right) & \text{if } p|l \text{ but } p \nmid l_1. \end{cases}$$

(*Note that $\eta(s;l)$ is absolutely convergent in $\mathrm{Re}(s) > 1/2$.*)

*Proof.* This follows by comparing the Euler factors on both sides.

From Lemma 5.1, we see that

(5.4a) $$\mathcal{P}_1(l) = \frac{2}{3}\frac{1+O(Y^{-1})}{\zeta(2)\sqrt{l_1}}d(l_1)I(l)$$

where

(5.4b) $$I(l) := \frac{1}{2\pi i}\int_{(c)}\left(\frac{\Gamma(\frac{s}{2}+\frac{1}{4})}{\Gamma(\frac{1}{4})}\right)^2\left(\frac{8X}{l_1\pi}\right)^s\check{\Phi}(s)\zeta(1+2s)^3\eta(1+2s;l)\frac{ds}{s}.$$

We move the line of integration in (5.4b) to the $\mathrm{Re}(s) = -\frac{1}{4}+\varepsilon$ line. There is a pole of order 4 at $s=0$ and we shall evaluate the residue of this pole shortly. We now bound the integral on the $-\frac{1}{4}+\varepsilon$ line. From [6, p. 79] we know that on this line $|\zeta(1+2s)| \ll |s|$, and plainly $|\eta(1+2s;l)| \ll \prod_{p|l_1}(1+O(\frac{1}{\sqrt{p}}))\prod_{p\nmid l_1}(1+O(\frac{1}{p^{1+\varepsilon}})) \ll l_1^\varepsilon$. Hence the integral on the $-\frac{1}{4}+\varepsilon$ line is
$$\ll \frac{l_1^{\frac{1}{4}+\varepsilon}}{X^{\frac{1}{4}-\varepsilon}}\int_{(-\frac{1}{4}+\varepsilon)}|\check{\Phi}(s)||s|^2|\Gamma(\tfrac{s}{2}+\tfrac{1}{4})|^2|ds| \ll \frac{l_1^{\frac{1}{4}+\varepsilon}}{X^{\frac{1}{4}-\varepsilon}}.$$

We now evaluate the residue of the pole at $s = 0$. For some absolute constants $c_1, c_2, \ldots, d_1, d_2, \ldots$, we have the Laurent series expansions

$$\left(\frac{\Gamma(s/2+1/4)}{\Gamma(1/4)}\right)^2 \zeta(1+2s)^3 = \frac{1}{8s^3} + \frac{c_1}{s^2} + \frac{c_2}{s} + c_3 + c_4 s + \ldots;$$

$$\left(\frac{8X}{l_1\pi}\right)^s = 1 + s\log\left(\frac{8X}{l_1\pi}\right) + \frac{s^2}{2!}\log^2\left(\frac{8X}{l_1\pi}\right) + \frac{s^3}{6}\log^3\left(\frac{8X}{l_1\pi}\right) + \ldots;$$

$$\eta(1+2s;l) = \eta(1;l)\left(1 + d_1 s \frac{\eta'}{\eta}(1;l) + d_2 s^2 \frac{\eta''}{\eta}(1;l) + d_3 s^3 \frac{\eta'''}{\eta}(1;l) + \ldots\right);$$

and $\check{\Phi}(s) = \check{\Phi}(0) + s\check{\Phi}'(0) + \frac{s^2}{2}\check{\Phi}''(0) + \ldots$. It follows that the residue may be written as

$$\frac{\eta(1;l)}{48}\check{\Phi}(0)\left(\log^3\left(\frac{8X}{l_1\pi}\right) + P_0\left(\log\frac{8X}{l_1\pi}\right) + \frac{\eta'}{\eta}(1;l)P_1\left(\log\frac{8X}{l_1\pi}\right)\right.$$
$$\left. + \frac{\eta''}{\eta}(1;l)P_2\left(\log\frac{8X}{l_1\pi}\right) + P_3\frac{\eta'''}{\eta}(1;l)\right),$$

where $P_3$ is an absolute constant, and $P_0$, $P_1$, and $P_2$ are polynomials of degrees 2, 2, and 1 respectively. Their coefficients involve absolute constants, and linear combinations of the parameters $\check{\Phi}^{(i)}(0)/\check{\Phi}(0)$ for $i = 1$, 2, and 3.

From the multiplicative definition of $\eta(s;l)$ in Lemma 5.1, we may write $\eta(s;l) = F(s)G_l(s)H_{l_1}(s)$ where $F(s)$ is independent of $l$; and $G_l(s) = \prod_{p|l} g_p(s)$ and $H_{l_1}(s) = \prod_{p|l_1} h_p(s)$ for appropriate Euler factors $g_p$ and $h_p$. Differentiating this product $i$ times, we see that $\frac{\eta^{(i)}}{\eta}(s;l)$ may be expressed as $\sum_{j,k=0}^{i} c_{j,k} \frac{G_l^{(j)}}{G_l}(1)\frac{H_{l_1}^{(k)}}{H_{l_1}}(1)$ for some absolute constants $c_{j,k}$ (given easily in terms of derivatives of $F(s)$). It is easy to see that $\frac{G_l^{(j)}}{G_l}(1) = \sum_{m|l} \frac{\Lambda_j(m)}{m} D_{1,j}(m)$, and that $\frac{H_{l_1}^{(k)}}{H_{l_1}}(1) = \sum_{n|l_1} \frac{\Lambda_k(n)}{n} D_{2,k}(n)$ where $D_{1,j}(m) \ll_j 1$ and $D_{2,k}(n) \ll_k 1$. From these observations, we may recast the residue of $I(l)$ above as

$$\frac{\eta(1;l)}{48}\check{\Phi}(0)\left(\log^3\left(\frac{X}{l_1}\right) + \mathcal{O}_1(l)\right)$$

where

$$\mathcal{O}_1(l) = \sum_{j,k=0}^{3} \sum_{m|l} \sum_{n|l_1} \frac{\Lambda_j(m)}{m} \frac{\Lambda_k(n)}{n} P_{j,k}\left(\log\frac{X}{l_1}\right) D_0(m,n),$$

where $D_0(m,n) \ll 1$, and the $P_{j,k}$ are polynomials of degree $\leq 2$ whose coefficients involve absolute constants and a linear combination of $\check{\Phi}^{(i)}(0)/\check{\Phi}(0)$ ($i = 1, 2, 3$).

Using this evaluation of $I(l)$ in (5.4a) we conclude that

$$\mathcal{P}_1(l) = \frac{\hat{\Phi}(0)\eta(1;l)}{72\zeta(2)} \frac{d(l_1)}{\sqrt{l_1}} \left( \log^3\left(\frac{X}{l_1}\right) + \mathcal{O}_1(l) \right)$$
$$+ O\left( \frac{d(l_1)\log^3 X}{\sqrt{l_1}Y} + \frac{1}{(l_1 X)^{\frac{1}{4}-\varepsilon}} \right).$$

5.2. *Extracting the secondary principal term from $\mathcal{R}_0(l)$.* For all real numbers $\xi$ and complex numbers $w$ with $\text{Re}(w) > 0$ we define

(5.5) $$f(\xi, w) = \int_0^\infty \tilde{F}_t\left(\frac{\xi}{t}\right) t^{w-1} dt.$$

Since $|\tilde{F}_t(\frac{\xi}{t})| \leq 2|\tilde{F}_t(0)| \ll e^{-\frac{t}{20X}}$ by Lemma 2.1, clearly the integral in (5.5) converges for $\text{Re}(w) > 0$. We now collect together some properties of $f(\xi, w)$.

LEMMA 5.2. *If $\xi \neq 0$ then*

$$f(\xi, w) = |\xi|^w \check{\Phi}(w) \int_0^\infty \omega_2\left(\frac{|\xi|\pi}{8Xz}\right) (\cos(2\pi z) + \text{sgn}(\xi)\sin(2\pi z)) \frac{dz}{z^{w+1}}.$$

*The integral above may be expressed as*

(5.6) $$\frac{1}{2\pi i} \int_{(c)} \left(\frac{\Gamma(\frac{s}{2} + \frac{1}{4})}{\Gamma(\frac{1}{4})}\right)^2 \left(\frac{8X}{|\xi|\pi}\right)^s (2\pi)^{-s+w} \Gamma(s-w)$$
$$\times \left( \cos(\tfrac{\pi}{2}(s-w)) + \text{sgn}(\xi)\sin(\tfrac{\pi}{2}(s-w)) \right) \frac{ds}{s},$$

*for any $\text{Re}(w) + 1 > c > \max(0, \text{Re}(w))$. For $\xi \neq 0$, $f(\xi, w)$ is a holomorphic function of $w$ in $\text{Re}(w) > -1$, and in the region $1 \geq \text{Re}(w) > -1$ satisfies the bound*

$$|f(\xi, w)| \ll (1 + |w|)^{-\text{Re}(w) - \frac{1}{2}} \exp\left( -\frac{1}{10} \frac{\sqrt{|\xi|}}{\sqrt{X(|w|+1)}} \right) |\xi|^w |\check{\Phi}(w)|.$$

*Proof.* From (5.5) and the definition of $\tilde{F}_t$ we have

$$f(\xi, w) = \int_0^\infty \left( \int_0^\infty F_t(y) \left( \cos(2\pi y \tfrac{\xi}{t}) + \sin(2\pi y \tfrac{\xi}{t}) \right) dy \right) t^{w-1} dt.$$

In the inner integral over $y$, we make the substitution $z = |\xi|y/t$, so that this integral becomes

$$\frac{t}{|\xi|} \int_0^\infty F_t\left(\frac{tz}{|\xi|}\right) (\cos(2\pi z) + \text{sgn}(\xi)\sin(2\pi z)) \, dz.$$

We use this above, and interchange the integrals over $z$ and $t$. Thus

$$f(\xi, w) = \frac{1}{|\xi|} \int_0^\infty \left( \int_0^\infty F_t\left(\frac{tz}{|\xi|}\right) t^w dt \right) (\cos(2\pi z) + \text{sgn}(\xi)\sin(2\pi z)) \, dz.$$

From the definition of $F_t$ the inner integral is

$$\omega_2\left(\frac{\pi|\xi|}{8Xz}\right) \int_0^\infty \Phi\left(\frac{tz}{|\xi|}\right) t^w dt = \omega_2\left(\frac{\pi|\xi|}{8Xz}\right)\left(\frac{|\xi|}{z}\right)^{w+1} \check{\Phi}(w),$$

by a change of variables. The first statement of the lemma follows at once.

By the definition of $\omega_2$ (see (1.5)) we have for any $c > 0$

$$\int_0^\infty \omega_2\left(\frac{|\xi|\pi}{8Xz}\right)(\cos(2\pi z) + \text{sgn}(\xi)\sin(2\pi z))\frac{dz}{z^{w+1}}$$
$$= \int_0^\infty (\cos(2\pi z) + \text{sgn}(\xi)\sin(2\pi z))\frac{1}{2\pi i}\int_{(c)}\left(\frac{\Gamma(\frac{s}{2}+\frac{1}{4})}{\Gamma(\frac{1}{4})}\right)^2\left(\frac{8X}{|\xi|\pi}\right)^s z^{s-w-1}\frac{ds}{s}dz.$$

If we choose $c$ so that $\text{Re}(w)+1 > c > \max(0, \text{Re}(w))$ (thus $0 < \text{Re}(s-w) < 1$) then we may interchange the two integrals above. This is because the $z$-integral $\int_0^\infty (\cos(2\pi z) + \text{sgn}(\xi)\sin(2\pi z))z^{s-w-1}dz$ is (conditionally) convergent for $\text{Re}(s)$ in this range. The interchange of integrals is rigourously justified by restricting the $z$-integral to the range $(\varepsilon, 1/\varepsilon)$ and letting $\varepsilon \to 0+$. Thus, with $c$ in this range, we have

$$\int_0^\infty \omega_2\left(\frac{|\xi|\pi}{8Xz}\right)(\cos(2\pi z) + \text{sgn}(\xi)\sin(2\pi z))\frac{dz}{z^{w+1}} = \frac{1}{2\pi i}\int_{(c)}\left(\frac{\Gamma(\frac{s}{2}+\frac{1}{4})}{\Gamma(\frac{1}{4})}\right)^2$$
$$\times \left(\frac{8X}{|\xi|\pi}\right)^s\left(\int_0^\infty (\cos(2\pi z) + \text{sgn}(\xi)\sin(2\pi z))z^{s-w-1}dz\right)\frac{ds}{s}.$$

Employing the expressions for the Fourier sine and cosine transforms of $z^{s-w-1}$ (see [7, pp. 1186–1190]), we obtain the second assertion of the lemma.

From the first two statements of the lemma, it is clear that for fixed $\xi \neq 0$, $f(\xi, w)$ is an analytic function of $w$ for $\text{Re}(w) > -1$; and it remains only to prove the bound on $|f(\xi, w)|$. Write the integral (5.6) as $\frac{1}{2\pi i}\int_{(c)} g(s, w; \text{sgn}(\xi))\left(\frac{8X}{|\xi|\pi}\right)^s ds$. Note that when $\text{Re}(w) > -1$, the integral (5.6) makes sense for $c > \max(0, \text{Re}(w))$. We shall bound it by choosing $c$ optimally. Stirling's formula shows that

$$|g(s, w; \text{sgn}(\xi))| \ll \left(\tfrac{|s|}{e}\right)^{c-\frac{3}{2}} \exp\left(-\tfrac{\pi}{2}|\text{Im}(s)|\right)(1 + |s-w|)^{c-\frac{1}{2}-\text{Re}(w)}.$$

Choosing $c = 1 + \frac{1}{3}\frac{\sqrt{|\xi|}}{\sqrt{X(1+|w|)}}$ we easily see the bound of the lemma.

Observe that for any sequence of numbers $a_n \ll n^\varepsilon$, and any smooth function $g$ with $g(0) = 0$ and $g(x)$ decaying rapidly as $x \to \infty$, we have the Mellin transform identity

$$\sum_{n=1}^\infty a_n g(n) = \frac{1}{2\pi i}\int_{(c)}\sum_{n=1}^\infty \frac{a_n}{n^w}\left(\int_0^\infty g(t)t^{w-1}dt\right)dw,$$

where $c > 1$. Hence we may recast the expression for $\mathcal{R}_0(l)$ (see (5.3) above) as

(5.7)
$$\mathcal{R}_0(l) = \frac{1}{2l} \sum_{\substack{\alpha \leq Y \\ (\alpha, 2l)=1}} \frac{\mu(\alpha)}{\alpha^2} \sum_{\substack{k=-\infty \\ k \neq 0}}^{\infty} \frac{(-1)^k}{2\pi i} \int_{(c)} \sum_{\substack{n=1 \\ (n,2\alpha)=1}}^{\infty} \frac{d(n)}{n^{\frac{3}{2}+w}} G_{4k}(ln) f\left(\frac{kX}{2\alpha^2 l}, w\right) dw,$$

for any $c > 0$.

LEMMA 5.3. *Write $4k = k_1 k_2^2$ where $k_1$ is a fundamental discriminant (possibly $k_1 = 1$ is the trivial character), and $k_2$ is positive. In the region $\mathrm{Re}(s) > 1$*

(5.8)
$$\sum_{\substack{n=1 \\ (n,2\alpha)=1}}^{\infty} \frac{d(n)}{n^s} \frac{G_{4k}(ln)}{\sqrt{n}} = L(s, \chi_{k_1})^2 \prod_p \mathcal{G}_p(s; k, l, \alpha) =: L(s, \chi_{k_1})^2 \mathcal{G}(s; k, l, \alpha)$$

*where $\mathcal{G}_p(s; k, l, \alpha)$ is defined as follows:*

$$\mathcal{G}_p(s; k, l, \alpha) = \left(1 - \frac{1}{p^s}\left(\frac{k_1}{p}\right)\right)^2 \quad \text{if } p|2\alpha, \quad \text{and}$$

$$\mathcal{G}_p(s; k, l, \alpha) = \left(1 - \frac{1}{p^s}\left(\frac{k_1}{p}\right)\right)^2 \sum_{r=0}^{\infty} \frac{d(p^r)}{p^{rs}} \frac{G_k(p^{r+\mathrm{ord}_p(l)})}{p^{\frac{r}{2}}}, \quad \text{if } p \nmid 2\alpha.$$

*Then $\mathcal{G}(s; k, l, \alpha)$ is holomorphic in the region $\mathrm{Re}(s) > \frac{1}{2}$, and for $\mathrm{Re}(s) \geq \frac{1}{2}+\varepsilon$ satisfies the bound*

(5.9)
$$|\mathcal{G}(s; k, l, \alpha)| \ll \alpha^\varepsilon |k|^\varepsilon l^{\frac{1}{2}+\varepsilon} (l, k_2^2)^{\frac{1}{2}}.$$

*Proof.* The Euler product expansion (5.8) follows from the multiplicativity of $G_{4k}(n)$ (see Lemma 2.3). By Lemma 2.3 we see that for a generic $p \nmid 2\alpha k l$, $\mathcal{G}_p(s; k, l, \alpha) = 1 - \frac{3}{p^{2s}} + \left(\frac{k_1}{p}\right)\frac{2}{p^{3s}}$. This shows that $\mathcal{G}(s; k, l, \alpha)$ is holomorphic in $\mathrm{Re}(s) > \frac{1}{2}$. It remains only to prove the bound (5.9). From our evaluation of $\mathcal{G}_p$ for $p \nmid 2kl\alpha$ we see that for $\mathrm{Re}(s) > \frac{1}{2} + \varepsilon$,

$$|\mathcal{G}(s; k, l, \alpha)| \ll (|k|l\alpha)^\varepsilon \prod_{\substack{p|kl \\ p \nmid 2\alpha}} |\mathcal{G}_p(s; k, l)|.$$

Suppose now that $p^a \parallel k$ and $p^b \parallel l$. Plainly we may suppose that $b \leq a+1$, else (using Lemma 2.3) $\mathcal{G}_p(s; k, l, \alpha) = 0$. Notice that $p^{[\frac{a}{2}]} \parallel k_2$. We now claim that $|\mathcal{G}_p(s; k, l, \alpha)| \ll (a+1)^2 p^{\min(b, [\frac{a}{2}]+\frac{b}{2})}$, which when inserted in our earlier estimate gives (5.9).

By the trivial bound $|G_k(p^r)| \leq p^r$ it follows that $|\mathcal{G}_p(s;k,l,\alpha)| \ll (a+1)^2 p^b$ so that our claim follows if $[\frac{a}{2}] \geq \frac{b}{2}$. The only remaining cases are $a$ even and $b = a+1$; and $a$ odd and $b = a$, or $b = a+1$. These are easily verified using Lemma 2.3.

We use Lemma 5.3 in (5.7), and move the line of integration to the line $\text{Re}(w) = -\frac{1}{2} + \varepsilon$. We encounter poles only when $k = \square$ (so that $k_1 = 1$, and $L(s, \chi_{k_1}) = \zeta(s)$); here there is a pole of order 2 at $w = 0$ and the residue is the source of the secondary principal term. Thus we may write $\mathcal{R}_0(l) = \mathcal{R}(l) + \mathcal{P}_2(l)$, where

$$(5.10) \quad \mathcal{R}(l) = \frac{1}{2l} \sum_{\substack{\alpha \leq Y \\ (\alpha, 2l)=1}} \frac{\mu(\alpha)}{\alpha^2} \sum_{\substack{k=-\infty \\ k \neq 0}}^{\infty} \frac{(-1)^k}{2\pi i}$$

$$\times \int_{(-\frac{1}{2}+\varepsilon)} L(1+w, \chi_{k_1})^2 \mathcal{G}(1+w; k, l, \alpha) f\left(\frac{kX}{2\alpha^2 l}, w\right) dw,$$

and (with an obvious change in notation, writing $k^2$ in place of $k$),

$$(5.11)$$
$$\mathcal{P}_2(l) = \frac{1}{2l} \operatorname*{Res}_{w=0} \sum_{\substack{\alpha \leq Y \\ (\alpha, 2l)=1}} \frac{\mu(\alpha)}{\alpha^2} \sum_{k=1}^{\infty} (-1)^k \zeta(1+w)^2 \mathcal{G}(1+w; k^2, l, \alpha) f\left(\frac{k^2 X}{2\alpha^2 l}, w\right).$$

5.3. *The secondary principal term $\mathcal{P}_2(l)$.* Below we shall suppose that $w$ is in the vicinity of 0; precisely, $|w| \leq \frac{1}{\log X}$. We begin by trying to simplify

$$(5.12) \quad \sum_{k=1}^{\infty} (-1)^k \mathcal{G}(1+w; k^2, l, \alpha) f\left(\frac{k^2 X}{2\alpha^2 l}, w\right).$$

We now define for $|v - 1| \leq \frac{1}{\log X}$, and any $u$ with $\text{Re}(u) > \frac{1}{2}$,

$$\mathcal{H}(u, v; l, \alpha) = l^u \sum_{k=1}^{\infty} \frac{(-1)^k}{k^{2u}} \mathcal{G}(v; k^2, l, \alpha).$$

Note that this series converges absolutely in this range. Define also

$$\Gamma_1(u) = (2\pi)^{-u} \Gamma(u) \left( \cos(\tfrac{\pi}{2} u) + \sin(\tfrac{\pi}{2} u) \right).$$

Using Lemma 5.2 and these definitions, we see that (5.12) may be recast as

$$\check{\Phi}(w) \left(\frac{X}{2\alpha^2}\right)^w \frac{1}{2\pi i} \int_{(\frac{3}{4})} \left(\frac{\Gamma(\frac{s}{2} + \frac{1}{4})}{\Gamma(\frac{1}{4})}\right)^2 \left(\frac{16\alpha^2}{\pi}\right)^s \Gamma_1(s-w) \mathcal{H}(s-w, 1+w; l, \alpha) \frac{ds}{s}.$$

From this it follows easily that

(5.13)
$$\operatorname*{Res}_{w=0} \zeta(1+w)^2 \times \text{(5.14)}$$

$$= \frac{\check{\Phi}(0)}{2\pi i} \int_{(\frac{3}{4})} \left(\frac{\Gamma(\frac{s}{2}+\frac{1}{4})}{\Gamma(\frac{1}{4})}\right)^2 \left(\frac{16\alpha^2}{\pi}\right)^s \Gamma_1(s)\mathcal{H}(s,1;l,\alpha)$$

$$\times \left(\log \frac{D_1 X}{\alpha^2} - \frac{\frac{\partial}{\partial s}\mathcal{H}(s,1;l,\alpha)}{\mathcal{H}(s,1;l,\alpha)} + \frac{\frac{\partial}{\partial w}\mathcal{H}(s,1+w;l,\alpha)}{\mathcal{H}(s,1;l,\alpha)}\bigg|_{w=0} - \frac{\Gamma'_1(s)}{\Gamma_1(s)}\right) \frac{ds}{s},$$

where $D_1$ is a constant depending only on $\Phi$. We note that $\log D_1$ may be written as $A + B\check{\Phi}'(0)/\check{\Phi}(0)$ for absolute constants $A$ and $B$.

From the definition of $\mathcal{G}(v; k^2, l, \alpha)$ we see that

$$\mathcal{H}(u,v;l,\alpha) = -l^u(1-2^{1-2u})\sum_{k=1}^\infty \frac{1}{k^{2u}}\mathcal{G}(v;k^2,l,\alpha)$$

$$= -l^u(1-2^{1-2u})\prod_p \sum_{b=0}^\infty \frac{\mathcal{G}_p(v;p^{2b},l,\alpha)}{p^{2bu}}.$$

Using the expression for $\mathcal{G}_p$ in Lemma 5.3, and then employing Lemma 2.3 to evaluate it, we may write

$$\mathcal{H}(u,v;l,\alpha) =: -l(1-2^{1-2u})l_1^{u-\frac{1}{2}}\zeta(2u)\zeta(2u+2v-1)\prod_p \mathcal{H}_{1,p}(u,v;l,\alpha)$$

$$=: -l(1-2^{1-2u})l_1^{u-\frac{1}{2}}\zeta(2u)\zeta(2u+2v-1)\mathcal{H}_1(u,v;l,\alpha),$$

where

$$\mathcal{H}_{1,p} = \begin{cases} \left(1-\frac{1}{p^v}\right)^2 \left(1-\frac{1}{p^{2u+2v-1}}\right) & \text{if } p|2\alpha \\ \frac{(1-\frac{1}{p^v})^2}{(1-\frac{1}{p^{2u+2v-1}})}\left(1+\frac{2}{p^v}-\frac{2}{p^{2u+v}}+\frac{1}{p^{2u+2v-1}}-\frac{3}{p^{2u+2v}}+\frac{1}{p^{4u+4v-1}}\right) & \text{if } p\nmid 2\alpha l \\ \frac{(1-\frac{1}{p^v})^2}{(1-\frac{1}{p^{2u+2v-1}})}\left(1-\frac{1}{p^{2u}}+\frac{2}{p^{2u+v-1}}-\frac{2}{p^{2u+v}}+\frac{1}{p^{2u+2v-1}}-\frac{1}{p^{4u+2v-1}}\right) & \text{if } p|l_1 \\ \frac{(1-\frac{1}{p^v})^2}{(1-\frac{1}{p^{2u+2v-1}})}\left(1-\frac{1}{p}+\frac{2}{p^v}-\frac{2}{p^{2u+v}}+\frac{1}{p^{2u+2v-1}}-\frac{1}{p^{2u+2v}}\right) & \text{if } p|l, p\nmid l_1. \end{cases}$$

From this it follows that

$$\mathcal{H}(s,1;l,\alpha), \quad \frac{\partial}{\partial s}\mathcal{H}(s,1;l,\alpha), \quad \text{and} \quad \frac{\partial}{\partial w}\mathcal{H}(s,w;l,\alpha)\bigg|_{w=1}$$

are holomorphic when $\operatorname{Re}(s) > 0$. Since $|\zeta(s)|, |\zeta'(s)| \ll \max((1+|s|)^\varepsilon, (1+|s|)^{\frac{1-\operatorname{Re}(s)}{2}+\varepsilon})$ (see [20, pp. 95, 96] for the proof of this estimate for $\zeta(s)$;

the estimate for $\zeta'(s)$ follows by a similar convexity principle) we obtain that when $\text{Re}(s) \geq \frac{1}{\log X}$ they are bounded by

$$(5.14) \qquad \ll l^{1+\varepsilon} l_1^{\text{Re}(s)-\frac{1}{2}} (\alpha X)^\varepsilon (|s|+1).$$

Hence we may move the line of integration in (5.13) to the $\text{Re}(s) = \frac{1}{\log X}$ line. We now introduce the sum over $\alpha$ as well, and (since $\check{\Phi}(0) = \hat{\Phi}(0)$) arrive at

$$\mathcal{P}_2(l) = \frac{\hat{\Phi}(0)}{2l} \frac{1}{2\pi i} \int_{(\frac{1}{\log X})} \left(\frac{\Gamma(\frac{s}{2}+\frac{1}{4})}{\Gamma(\frac{1}{4})}\right)^2 \left(\frac{16}{\pi}\right)^s \Gamma_1(s) \sum_{\substack{\alpha \leq Y \\ (\alpha, 2l)=1}} \frac{\mu(\alpha)}{\alpha^{2-2s}} \mathcal{H}(s, 1; l, \alpha)$$

$$\times \left(\log\left(\frac{D_1 X}{\alpha^2}\right) - \frac{\frac{\partial}{\partial s}\mathcal{H}(s,1;l,\alpha)}{\mathcal{H}(s,1;l,\alpha)} + \frac{\frac{\partial}{\partial w}\mathcal{H}(s,w;l,\alpha)}{\mathcal{H}(s,w;l,\alpha)}\bigg|_{w=1} - \frac{\Gamma_1'(s)}{\Gamma_1(s)}\right) \frac{ds}{s}.$$

We now extend the sum over $\alpha$ above to infinity. By (5.14) the error incurred in doing so is

$$\ll \frac{l^\varepsilon}{\sqrt{l_1} Y^{1-\varepsilon}} \int_{(\frac{1}{\log X})} |\Gamma(\tfrac{s}{2}+\tfrac{1}{4})|^2 \max(|\Gamma_1(s)|, |\Gamma_1'(s)|)|ds| \ll \frac{l^\varepsilon}{\sqrt{l_1} Y^{1-\varepsilon}},$$

because the integrand decays exponentially as $|\text{Im}(s)| \to \infty$. Hence

$$\mathcal{P}_2(l) = O\left(\frac{l^\varepsilon}{\sqrt{l_1} Y^{1-\varepsilon}}\right) + \frac{\hat{\Phi}(0)}{2\pi i} \int_{(\frac{1}{\log X})} \left(\frac{\Gamma(\frac{s}{2}+\frac{1}{4})}{\Gamma(\frac{1}{4})}\right)^2 \left(\frac{16}{\pi}\right)^s \Gamma_1(s) \mathcal{K}(s,1;l)$$

$$\times \left(\log(D_1 X) - \frac{\frac{\partial}{\partial s}\mathcal{K}(s,1;l)}{\mathcal{K}(s,1;l)} + \frac{\frac{\partial}{\partial w}\mathcal{K}(s,w;l)}{\mathcal{K}(s,w;l)}\bigg|_{w=1} - \frac{\Gamma_1'(s)}{\Gamma_1(s)}\right) \frac{ds}{s},$$

where

$$\mathcal{K}(s,w;l) = \frac{1}{2l} \sum_{\substack{\alpha=1 \\ (\alpha,2l)=1}}^{\infty} \frac{\mu(\alpha)}{\alpha^{2-2s}} \mathcal{H}(s,w;l,\alpha).$$

By our expression for $\mathcal{H}$, a calculation gives

$$\mathcal{K}(s,1;l) = -\frac{1}{8\sqrt{l_1}} \frac{\varphi(l)^2}{l^2} \prod_{\substack{p|l \\ p \nmid l_1}} \left(1+\frac{1}{p}\right) \frac{4^s + 4^{-s} - \frac{5}{2}}{4^s} \sum_{ab=l_1} \left(\frac{a}{b}\right)^s \zeta(2s)\zeta(2s+1)$$

$$\times \prod_{p \nmid 2l} \left(1-\frac{1}{p}\right)^2 \left(1+\frac{2}{p}+\frac{1}{p^3}-\frac{1}{p^2}(p^{-2s}+p^{2s})\right).$$

Using this together with the functional equation for $\zeta(s)$ and the relations $\Gamma(z)\Gamma(1-z) = \pi\text{cosec}(\pi z)$ and $\Gamma(z)\Gamma(z+\frac{1}{2}) = \pi^{\frac{1}{2}} 2^{1-2z} \Gamma(2z)$ we see that

(5.15)
$$\left(\frac{\Gamma(\frac{s}{2}+\frac{1}{4})}{\Gamma(\frac{1}{4})}\right)^2 \left(\frac{16}{\pi}\right)^s \Gamma_1(s) \mathcal{K}(s,1;l) = \left(\frac{\Gamma(-\frac{s}{2}+\frac{1}{4})}{\Gamma(\frac{1}{4})}\right)^2 \left(\frac{16}{\pi}\right)^{-s} \Gamma_1(-s) \mathcal{K}(-s,1;l).$$

Further, after more calculations, we have

$$\frac{\frac{\partial}{\partial w}\mathcal{K}(s,w;l)}{\mathcal{K}(s,w;l)}\bigg|_{w=1} - \frac{\frac{\partial}{\partial s}\mathcal{K}(s,1;l)}{\mathcal{K}(s,1;l)} = -\log l_1 - 2\frac{\zeta'}{\zeta}(2s) + 2\frac{\zeta'}{\zeta}(2s+1) + \Psi(s)$$

where

$$\Psi(s) = 2\log 2 + \frac{6\log 2}{(1-2^{1+2s})(1-2^{1-2s})} + \sum_{p|l}\frac{2\log p}{p-1} - \sum_{\substack{p|l \\ p\nmid l_1}}\frac{2\log p}{p+1}$$

$$+ \sum_{p\nmid 2l}\left(\frac{2\log p}{p-1} - \frac{2\log p}{p}\frac{1+\frac{2}{p^2}-\frac{1}{p}(p^{-2s}+p^{2s})}{1+\frac{2}{p}+\frac{1}{p^3}-\frac{1}{p^2}(p^{2s}+p^{-2s})}\right)$$

$$= \Psi(-s).$$

Note also that

$$\frac{\Gamma_1'(s)}{\Gamma_1(s)} = -\log(2\pi) + \frac{\Gamma'}{\Gamma}(s) + \frac{\pi}{2}\frac{\cos(\frac{\pi s}{2}) - \sin(\frac{\pi s}{2})}{\cos(\frac{\pi s}{2}) + \sin(\frac{\pi s}{2})}.$$

Using these together with the logarithmic derivative of the functional equation for $\zeta(s)$ and the logarithmic derivative of the relation $\Gamma(z)\Gamma(1-z) = \pi\operatorname{cosec}(\pi z)$ we conclude that

$$(5.16) \quad -\frac{\frac{\partial}{\partial s}\mathcal{K}(s,1;l)}{\mathcal{K}(s,1;l)} + \frac{\frac{\partial}{\partial w}\mathcal{K}(s,w;l)}{\mathcal{K}(s,w;l)}\bigg|_{w=1} - \frac{\Gamma_1'(s)}{\Gamma_1(s)}$$

$$= -\frac{\frac{\partial}{\partial s}\mathcal{K}(-s,1;l)}{\mathcal{K}(-s,1;l)} + \frac{\frac{\partial}{\partial w}\mathcal{K}(-s,w;l)}{\mathcal{K}(-s,w;l)}\bigg|_{w=1} - \frac{\Gamma_1'(-s)}{\Gamma_1(-s)}.$$

The evenness of the expressions in (5.15) and (5.16) is extremely convenient below. However we do not know of any 'natural' proofs of these facts, which perhaps make them look purely a matter of good fortune.

We now use these explicit calculations to complete our evaluation of $\mathcal{P}_2(l)$. Write

$$\mathcal{J}(s,l) = \left(\frac{\Gamma(\frac{s}{2}+\frac{1}{4})}{\Gamma(\frac{1}{4})}\right)^2 \left(\frac{16}{\pi}\right)^s \Gamma_1(s)\mathcal{K}(s,1;l)$$

$$\times \left(\log(D_1 X) - \frac{\frac{\partial}{\partial s}\mathcal{K}(s,1;l)}{\mathcal{K}(s,1;l)} + \frac{\frac{\partial}{\partial w}\mathcal{K}(s,w;l)}{\mathcal{K}(s,w;l)}\bigg|_{w=1} - \frac{\Gamma_1'(s)}{\Gamma_1(s)}\right),$$

so that by (5.15) and (5.16), $\mathcal{J}(s,l) = \mathcal{J}(-s,l)$. Now

$$\mathcal{P}_2(l) = O\left(\frac{l^\varepsilon}{\sqrt{l_1}Y^{1-\varepsilon}}\right) + \frac{\hat{\Phi}(0)}{2\pi i}\int_{(\frac{1}{\log X})}\mathcal{J}(s,l)\frac{ds}{s}.$$

We move the line of integration to the line $-\frac{1}{\log X} - i\infty$ to $-\frac{1}{\log X} + i\infty$, encountering a pole at $s=0$. Thus

$$\mathcal{P}_2(l) = O\left(\frac{l^\varepsilon}{\sqrt{l_1}Y^{1-\varepsilon}}\right) + \hat{\Phi}(0)\operatorname*{Res}_{s=0}\frac{\mathcal{J}(s,l)}{s} + \frac{\hat{\Phi}(0)}{2\pi i}\int_{-(\frac{1}{\log X})}\mathcal{J}(s,l)\frac{ds}{s},$$

and changing $s$ to $-s$ and using the relation $\mathcal{J}(s,l) = \mathcal{J}(-s,l)$ we see that the above is

$$= O\left(\frac{l^\varepsilon}{\sqrt{l_1}Y^{1-\varepsilon}}\right) + \hat{\Phi}(0)\operatorname*{Res}_{s=0}\frac{\mathcal{J}(s,l)}{s} - \frac{\hat{\Phi}(0)}{2\pi i}\int_{(\frac{1}{\log X})}\mathcal{J}(s,l)\frac{ds}{s}.$$

Hence

$$\mathcal{P}_2(l) = O\left(\frac{l^\varepsilon}{\sqrt{l_1}Y^{1-\varepsilon}}\right) + \frac{\hat{\Phi}(0)}{2}\operatorname*{Res}_{s=0}\frac{\mathcal{J}(s,l)}{s}.$$

To compute the residue above we shall employ the following Laurent series expansions. In these expansions, we shall use the symbol $O(s^n)$ to group together terms involving at least the $n$-th power of $s$. Note that

$$\left(\frac{\Gamma(\frac{s}{2}+\frac{1}{4})}{\Gamma(\frac{1}{4})}\right)^2\left(\frac{16}{\pi}\right)^s\Gamma_1(s)\frac{4^s+4^{-s}-\frac{5}{2}}{4^s}\zeta(2s)\zeta(2s+1) = \frac{1}{8s^2} + a_0 + O(s),$$

for an absolute constant $a_0$. Further, since $l_1$ is square-free,

$$\sum_{ab=l_1}\left(\frac{a}{b}\right)^s = \prod_{p|l_1}(p^{-s}+p^s) = \prod_{p|l_1}(2 + s^2\log^2 p + O(s^4))$$

$$= d(l_1)\left(1 + \frac{s^2}{2}\sum_{p|l_1}\log^2 p + O(s^4)\right).$$

Next

$$\prod_{p\nmid 2l}\left(1-\frac{1}{p}\right)^2\left(1+\frac{2}{p}+\frac{1}{p^3}-\frac{1}{p^2}(p^{-2s}+p^{2s})\right)$$

$$= \prod_{p\nmid 2l}\left(1-\frac{1}{p}\right)^2\left(1+\frac{2}{p}+\frac{1}{p^3}-\frac{2}{p^2}\right)\left(1 + s^2\sum_{p\nmid 2l}\frac{\log^2 p}{p^2}D(p) + O(s^4)\right)$$

for some $D(p) \ll 1$; and this may be rewritten as

$$\prod_{p\nmid 2l}\left(1-\frac{1}{p}\right)^2\left(1+\frac{2}{p}+\frac{1}{p^3}-\frac{2}{p^2}\right)\left(1 + s^2\left(a_1 - \sum_{p|l}\frac{\log^2 p}{p^2}D(p)\right) + O(s^4)\right),$$

for some absolute constant $a_1$. Lastly we note that in view of (5.16)

$$\left(\log(D_1 X) - \frac{\frac{\partial}{\partial s}\mathcal{K}(s,1;l)}{\mathcal{K}(s,1;l)} + \frac{\frac{\partial}{\partial w}\mathcal{K}(s,w;l)}{\mathcal{K}(s,w;l)}\bigg|_{w=1} - \frac{\Gamma_1'(s)}{\Gamma_1(s)}\right)$$

is an even function of $s$, and so its Laurent expansion involves only even powers of $s$. A little calculation shows that this expansion may be written as

$$\log\left(\frac{D_1 X}{l_1}\right) + a_2 + \Psi(0) + \frac{s^2}{2}(\Psi''(0) + a_3) + O(s^4),$$

for some absolute constants $a_2$ and $a_3$.

From these observations we see that for absolute constants $a_4, a_5, \ldots$,

$$\operatorname*{Res}_{s=0} \frac{\mathcal{J}(s,l)}{s} = -\frac{d(l_1)}{128\sqrt{l_1}} \frac{\varphi(l)^2}{l^2} \prod_{\substack{p|l \\ p \nmid l_1}} \left(1 + \frac{1}{p}\right) \prod_{p \nmid 2l} \left(1 - \frac{1}{p}\right)^2 \left(1 + \frac{2}{p} - \frac{2}{p^2} + \frac{1}{p^3}\right)$$

$$\times \left(\log \frac{D_1 X}{l_1} \sum_{p|l_1} \log^2 p + a_4 \log \frac{D_1 X}{l_1} + a_5 \right.$$

$$\left. + a_6 \Psi(0) + a_7 \Psi''(0) + a_8 \sum_{p|l} \frac{\log^2 p}{p^2} D(p) \right).$$

From its definition we see that $\Psi(0) = D_0 + \sum_{p|l_1} \frac{\log p}{p} D_2(p) + \sum_{p|l} \frac{\log p}{p} D_3(p)$ for an absolute constant $D_0$ and with $D_2(p), D_3(p) \ll 1$. Further, we may write $\Psi''(0)$ as $D_4 + \sum_{p|l} \frac{\log p}{p} D_5(p)$ for an absolute constant $D_4$, and with $D_5(p) \ll 1$. From these remarks, and keeping in mind the definition of $\eta(1; l)$, we get that

$$\operatorname*{Res}_{s=0} \frac{\mathcal{J}(s,l)}{s} = -\frac{d(l_1)}{\sqrt{l_1}} \frac{\eta(1;l)}{12\zeta(2)} \left(\log \frac{X}{l_1} \sum_{p|l_1} \log^2 p + \mathcal{O}_2(l) \right),$$

where

$$\mathcal{O}_2(l) = P\left(\log \frac{X}{l_1}\right) + \log D_1 \sum_{p|l_1} \log^2 p + \sum_{p|l} \frac{\log p}{p} D_6(p) + \sum_{p|l_1} \frac{\log p}{p} D_7(p)$$

with $P$ a polynomial of degree 1 whose coefficients involve absolute constants and the parameter $\check{\Phi}'(0)/\check{\Phi}(0)$; and $D_6(p)$ and $D_7(p)$ are $\ll 1$. We also recall here that $\log D_1$ may be written as $A + B\check{\Phi}'(0)/\check{\Phi}(0)$ for absolute constants $A$ and $B$.

We conclude from the above calculations that

$$\mathcal{P}_2(l) = -\hat{\Phi}(0) \frac{d(l_1)}{\sqrt{l_1}} \frac{\eta(1;l)}{24\zeta(2)} \left(\log \frac{X}{l_1} \sum_{p|l_1} \log^2 p + \mathcal{O}_2(l)\right) + O\left(\frac{l^\varepsilon}{\sqrt{l_1} Y^{1-\varepsilon}}\right).$$

5.4. *The contribution of the remainder terms $\mathcal{R}(l)$.* We begin by obtaining a bound for $\mathcal{R}(l)$ for individual $l$. Using the bounds of Lemmas 5.2 and 5.3 in (5.10) we obtain that $|\mathcal{R}(l)|$ is

$$\ll \frac{l^\varepsilon}{X^{\frac{1}{2}-\varepsilon}} \sum_{\alpha \leq Y} \frac{1}{\alpha^{1-\varepsilon}} \int_{(-\frac{1}{2}+\varepsilon)} |\check{\Phi}(w)|$$

$$\times \sum_{\substack{k=-\infty \\ k \neq 0}}^{\infty} \frac{|L(1+w, \chi_{k_1})|^2}{\sqrt{|k_1|}} \exp\left(-\frac{1}{10} \frac{\sqrt{|k|}}{\alpha\sqrt{l}(|w|+1)}\right) |dw|.$$

Performing the sum over $k_2$, we see that this is bounded by

$$\frac{l^{\frac{1}{2}+\varepsilon}}{X^{\frac{1}{2}-\varepsilon}} \sum_{\alpha \leq Y} \alpha^\varepsilon \int_{(-\frac{1}{2}+\varepsilon)} |\check{\Phi}(w)|(|w|+1)^{\frac{1}{2}}$$
$$\times \sum_{k_1} \frac{|L(1+w,\chi_{k_1})|^2}{|k_1|} \exp\left(-\frac{1}{10}\frac{\sqrt{|k_1|}}{\alpha\sqrt{l(|w|+1)}}\right)|dw|.$$

Splitting the $k_1$ into dyadic blocks, and using Lemma 2.5 to estimate these contributions we deduce that

$$|\mathcal{R}(l)| \ll \frac{l^{\frac{1}{2}+\varepsilon}Y^{1+\varepsilon}}{X^{\frac{1}{2}-\varepsilon}} \int_{(-\frac{1}{2}+\varepsilon)} |\check{\Phi}(w)|(|w|+1)^{1+\varepsilon}|dw| \ll \frac{l^{\frac{1}{2}+\varepsilon}Y^{1+\varepsilon}}{X^{\frac{1}{2}-\varepsilon}}\Phi_{(2)}\Phi_{(3)}^\varepsilon,$$

where the last inequality follows by using (1.4) with $\nu = 2$ for $|w| \leq \frac{\Phi_{(3)}}{\Phi_{(2)}}$, and $\nu = 3$ for larger $|w|$.

We now show how a better bound for $\mathcal{R}(l)$ may be obtained on average. Let $\beta_l = \frac{\overline{\mathcal{R}(l)}}{|\mathcal{R}(l)|}$ if $\mathcal{R}(l) \neq 0$, and $\beta_l = 1$ otherwise. Then, from (5.10),

(5.17)
$$\sum_{l=L}^{2L-1} |\mathcal{R}(l)| = \sum_{l=L}^{2L-1} \beta_l \mathcal{R}(l) \ll \sum_{\substack{\alpha \leq Y \\ (\alpha,2)=1}} \frac{1}{\alpha^2} \sum_{\substack{k=-\infty \\ k \neq 0}}^{\infty} \int_{(-\frac{1}{2}+\varepsilon)} |L(1+w,\chi_{k_1})|^2$$
$$\times \Bigg| \sum_{\substack{l=L \\ (l,\alpha)=1}}^{2L-1} \frac{\beta_l}{l} \mathcal{G}(1+w; k, l, \alpha) f\left(\frac{kX}{2\alpha^2 l}, w\right) \Bigg| |dw|.$$

We now split the sum over $k$ into dyadic blocks $K \leq |k| \leq 2K-1$. By Cauchy's inequality

$$\sum_{|k|=K}^{2K-1} |L(1+w,\chi_{k_1})|^2 \Bigg| \sum_{\substack{l=L \\ (l,\alpha)=1}}^{2L-1} \frac{\beta_l}{l} \mathcal{G}(1+w; k, l, \alpha) f\left(\frac{kX}{2\alpha^2 l}, w\right) \Bigg|$$
$$\ll \left( \sum_{|k|=K}^{2K-1} k_2 |L(1+w,\chi_{k_1})|^4 \right)^{\frac{1}{2}}$$
$$\times \left( \sum_{|k|=K}^{2K-1} \frac{1}{k_2} \Bigg| \sum_{\substack{l=L \\ (l,\alpha)=1}}^{2L-1} \frac{\beta_l}{l} \mathcal{G}(1+w; k, l, \alpha) f\left(\frac{kX}{2\alpha^2 l}, w\right) \Bigg|^2 \right)^{\frac{1}{2}},$$

and using Lemma 2.5 to estimate the first factor, this is

(5.18)
$$\ll (K(1+|w|))^{\frac{1}{2}+\varepsilon} \left( \sum_{|k|=K}^{2K-1} \frac{1}{k_2} \Bigg| \sum_{\substack{l=L \\ (l,\alpha)=1}}^{2L-1} \frac{\beta_l}{l} \mathcal{G}(1+w; k, l, \alpha) f\left(\frac{kX}{2\alpha^2 l}, w\right) \Bigg|^2 \right)^{\frac{1}{2}}.$$

LEMMA 5.4. *Let $\alpha \leq Y$, $K$ and $L$ be positive integers, and suppose $w$ is a complex number with $\operatorname{Re}(w) = -\frac{1}{2} + \varepsilon$. Then for any choice of complex numbers $\gamma_l$ with $|\gamma_l| \leq 1$,*

$$\sum_{|k|=K}^{2K-1} \frac{1}{k_2} \left| \sum_{l=L}^{2L-1} \frac{\gamma_l}{l} \mathcal{G}(1+w; k, l, \alpha) f\left(\frac{kX}{2\alpha^2 l}, w\right) \right|^2$$

*is bounded by*

$$(1+|w|)^\varepsilon |\check{\Phi}(w)|^2 \frac{\alpha^{2+\varepsilon} L^{2+\varepsilon} K^\varepsilon}{X^{1-\varepsilon}} \exp\left(-\frac{1}{20} \frac{\sqrt{K}}{\alpha \sqrt{L(1+|w|)}}\right),$$

*and also by*

$$((1+|w|)\alpha KLX)^\varepsilon |\check{\Phi}(w)|^2 \frac{\alpha^2(KL + L^2)}{KX}.$$

Before proving this lemma we note the bound it gives for $\sum_{l=L}^{2L-1} |\mathcal{R}(l)|$. We bound (5.18) using the first bound of the lemma for $K \geq \alpha^2 L(1+|w|) \log^2 X$, and the second bound for smaller $K$. Inserting this bound in (5.17) gives (with a little calculation)

$$\sum_{l=L}^{2L-1} |\mathcal{R}(l)| \ll \frac{L^{1+\varepsilon}}{X^{\frac{1}{2}-\varepsilon}} \sum_{\alpha \leq Y} \alpha^\varepsilon \int_{(-\frac{1}{2}+\varepsilon)} |\check{\Phi}(w)|(1+|w|)^{1+\varepsilon} |dw|$$

$$\ll \frac{L^{1+\varepsilon} Y^{1+\varepsilon}}{X^{\frac{1}{2}-\varepsilon}} \Phi_{(2)} \Phi_{(3)}^\varepsilon,$$

as desired.

*Proof.* Using the bound for $\mathcal{G}$ in Lemma 5.3, and the bound for $|f(\xi, w)|$ in Lemma 5.2 we obtain that our desired sum is

$$\ll (1+|w|)^\varepsilon |\check{\Phi}(w)|^2 \exp\left(-\frac{1}{20} \frac{\sqrt{K}}{\alpha \sqrt{L(1+|w|)}}\right)$$

$$\times \frac{\alpha^{2+\varepsilon}(LK)^\varepsilon}{X^{1-\varepsilon}} \sum_{|k|=K}^{2K-1} \frac{1}{|k|k_2} \left| \sum_{l=L}^{2L} (l, k_2^2)^{\frac{1}{2}} \right|^2.$$

This immediately gives the first bound of the lemma.

Write the integral (5.6) as $\frac{1}{2\pi i} \int_{(c)} g(s, w; \operatorname{sgn}(\xi)) \left(\frac{8X}{|\xi|\pi}\right)^s ds$. Taking $c = \varepsilon$, we see that

$$\left| \sum_{l=L}^{2L-1} \frac{\gamma_l}{l} \mathcal{G}(1+w; k, l, \alpha) f\left(\frac{kX}{2\alpha^2 l}, w\right) \right|$$

$$\ll |\check{\Phi}(w)| \frac{\alpha^{1+\varepsilon}}{(|k|X)^{\frac{1}{2}-\varepsilon}} \int_{(c)} \left| g(s, w; \operatorname{sgn}(k)) \sum_{l=L}^{2L-1} \frac{\gamma_l}{l^{1+w-s}} \mathcal{G}(1+w; k, l, \alpha) ds \right|.$$

Since $|g(s, w; \text{sgn}(k))| \ll (1 + |w|)^\varepsilon \exp(-\frac{\pi}{2}|\text{Im}(s)|)$ by Stirling's formula, we get by Cauchy's inequality that

$$\Big|\sum_{l=L}^{2L-1} \frac{\gamma_l}{l} \mathcal{G}(1+w; k, l, \alpha) f\Big(\frac{kX}{2\alpha^2 l}, w\Big)\Big|^2 \ll (1+|w|)^\varepsilon |\check{\Phi}(w)|^2 \frac{\alpha^{2+\varepsilon}}{(|k|X)^{1-\varepsilon}}$$

$$\times \int_{(c)} \exp\left(-\frac{\pi}{2}|\text{Im}(s)|\right) \Big|\sum_{l=L}^{2L-1} \frac{\gamma_l}{l^{1+w-s}} \mathcal{G}(1+w; k, l, \alpha)\Big|^2 |ds|.$$

The second bound of the lemma follows by combination of this with Lemma 5.5 below.

LEMMA 5.5. *Let $|\delta_l| \ll l^\varepsilon$ be any sequence of complex numbers and let $w$ be any complex number with $\text{Re}(w) = -\frac{1}{2} + \varepsilon$. Then*

$$\sum_{|k|=K}^{2K-1} \frac{1}{k_2} \Big|\sum_{l=L}^{2L-1} \frac{\delta_l}{\sqrt{l}} \mathcal{G}(1+w; k, l, \alpha)\Big|^2 \ll (KL\alpha)^\varepsilon (K+L)L.$$

*Proof.* For any integer $k = \pm \prod_{i, a_i \geq 1} p_i^{a_i}$ we define $a(k) = \prod_i p_i^{a_i+1}$, and put $b(k) = \prod_{i, a_i=1} p_i \prod_{i, a_i \geq 2} p_i^{a_i-1}$. Note that $\mathcal{G}(1+w; k, l, \alpha) = 0$ unless $l$ can be written as $dm$ where $d|a(k)$ and $(m, k) = 1$ with $m$ square-free. From the definition of $\mathcal{G}$ in Lemma 5.3, and using Lemma 2.3, we get

$$\mathcal{G}(1+w; k, l, \alpha) = \sqrt{m} \Big(\frac{k}{m}\Big) \prod_{p|m} \Big(1 + \frac{2}{p^{1+w}} \Big(\frac{k}{p}\Big)\Big)^{-1} \mathcal{G}(1+w; k, d, \alpha).$$

Using Lemma 5.3 to bound $|\mathcal{G}(1+w; k, d, \alpha)|$ we see that our desired sum is

$$\ll (KL\alpha)^\varepsilon \sum_{|k|=K}^{2K-1} \frac{1}{k_2} \sum_{d|a(k)} d \Big|\sum_{m=L/d}^{2L/d} \delta_{dm} \mu(m)^2 \Big(\frac{k}{m}\Big) \prod_{p|m} \Big(1 + \frac{2}{p^{1+w}} \Big(\frac{k}{p}\Big)\Big)^{-1}\Big|^2.$$

We interchange the sums over $d$ and $k$. Note that $d|a(k)$ implies that $b(d)|k$, so that $k = b(d)f$ for some integer $f$ with $K/b(d) \leq |f| \leq 2K/b(d)$. Write $4f = f_1 f_2^2$ where $f_1$ is a fundamental discriminant, and $f_2$ is positive. Notice that $k_2 \geq f_2$. Thus our desired sum is bounded by

$$(KL\alpha)^\varepsilon \sum_{d \leq 2L} d \sum_{f=K/b(d)}^{2K/b(d)} \frac{1}{f_2} \Big|\sum_{m=L/d}^{2L/d} \delta(dm) \mu(m)^2 \Big(\frac{fb(d)}{m}\Big) \prod_{p|m} \Big(1 + \frac{2}{p^{1+w}} \Big(\frac{fb(d)}{p}\Big)\Big)^{-1}\Big|^2,$$

and by Lemma 2.4 this is

$$\ll (KL\alpha)^\varepsilon \sum_{d \leq L} d \frac{L}{d} \Big(\frac{K}{b(d)} + \frac{L}{d}\Big) \ll (KL\alpha)^\varepsilon (KL + L^2) \sum_{d \leq L} \frac{1}{b(d)}$$

$$\ll (KL\alpha)^\varepsilon (KL + L^2).$$

5.5. *Proofs of Proposition* 1.3 *and Corollary* 1.4. We assemble the asymptotic formulae for $\mathcal{P}_1(l)$ and $\mathcal{P}_2(l)$ of Section 5.1, and Section 5.3, and the bounds on $\mathcal{R}(l)$ in Section 5.4, and check that $\eta(1;l) = Dl_1/(\sigma(l_1)h(l))$. This proves Proposition 1.3.

Notice that by Lemma 2.2, Proposition 1.1 with $M(d) = 1$, and Proposition 1.3 (with $l = 1$),

$$\mathcal{S}(L(\tfrac{1}{2},\chi_{8d})^2;\Phi) = \hat{\Phi}(0)Q_{0;\Phi}(\log X) + O\left(\frac{X^\varepsilon}{Y} + \Phi_{(2)}\Phi_{(3)}^\varepsilon \frac{Y^{1+\varepsilon}}{X^{\frac{1}{2}-\varepsilon}}\right),$$

where $Q_{0;\Phi}$ is a polynomial of degree 3 whose coefficients involve the parameters $\frac{\check{\Phi}^{(j)}(0)}{\check{\Phi}(0)}$ for $j = 1,2,3$. Now we choose $\Phi$ such that $\Phi(t) = 1$ for $t \in (1+Z^{-1}, 2-Z^{-1})$ and such that $\Phi^{(\nu)}(t) \ll_\nu Z^\nu$, for all $\nu \geq 0$. It follows that $\Phi_{(2)} \ll Z$, $\Phi_{(3)} \ll Z^2$, and that $\check{\Phi}(0) = \hat{\Phi}(0) = 1 + O(Z^{-1})$. Further, the parameters $\check{\Phi}^{(j)}(0)/\check{\Phi}(0)$ equal $\int_1^2 (\log y)^j dy + O(Z^{-1})$. Thus we deduce that for a polynomial $Q_0$ whose coefficients are absolute constants

$$\sum_{X \leq d \leq 2X} \Phi(\tfrac{d}{X})L(\tfrac{1}{2},\chi_{8d})^2 = XQ_0(\log X) + O\left(\frac{X^{1+\varepsilon}}{Y} + \frac{X^{1+\varepsilon}}{Z} + Z^{1+\varepsilon}Y^{1+\varepsilon}X^{\frac{1}{2}+\varepsilon}\right).$$

We now apply Lemma 2.2, Propositions 1.1 and 1.3 to the new choice $\Phi_1(t) = 1 - \Phi(t)$ for $1 \leq t \leq 2$ and $\Phi_1(t) = 0$ otherwise. Then we see that

$$\sum_{X \leq d \leq 2X} (1 - \Phi(\tfrac{d}{X}))L(\tfrac{1}{2},\chi_{8d})^2 \ll \frac{X^{1+\varepsilon}}{Y} + \frac{X^{1+\varepsilon}}{Z} + Y^{1+\varepsilon}Z^{1+\varepsilon}X^{\frac{1}{2}+\varepsilon}.$$

We add the two displays above, and take $Y = Z = X^{\frac{1}{6}}$ to obtain

$$\sum_{X \leq d \leq 2X} L(\tfrac{1}{2},\chi_{8d})^2 = XQ_0(\log X) + O(X^{\frac{5}{6}+\varepsilon}).$$

Having summed this with $X = x/2$, $X = x/4$, ..., we have proved Corollary 1.4.

## 6. Choosing the mollifier: proof of Theorem 1

Throughout we shall suppose that $\lambda(l) \ll l^{-1+\varepsilon}$ and that $Y = X^\varepsilon$. We choose $M = (\sqrt{X})^\theta$ where $\theta \leq 1 - \varepsilon$. For simplicity we shall suppose that $\lambda(l) = 0$ unless $l$ is odd and square-free. The optimal mollifier satisfies these constraints, so no loss of generality is incurred in making this simplification. Further, put $\lambda(l) = 0$ if $l > M$. Lastly, we shall suppose that $\Phi$ is chosen so that $\Phi_{(j)} \ll_j X^\varepsilon$ for all $j \geq 1$; indeed later we shall just choose $\Phi$ to be any smooth approximation to the characteristic function of $(1,2)$.

With hindsight, we introduce the linear change of variables

$$\xi(\gamma) = \sum_a \frac{\lambda(a\gamma)}{h(a)} \frac{ad(a)}{\sigma(a)}. \tag{6.1}$$

Like $\lambda$, $\xi$ is supported only on odd square-free integers below $M$. This change of variables is invertible, and $\lambda$ may be recovered from $\xi$ by

$$\lambda(l) = \sum_a \frac{\mu(a)}{h(a)} \frac{ad(a)}{\sigma(a)} \xi(la). \tag{6.2}$$

We shall further require our mollifier to satisfy

$$|\xi(\gamma)| \ll \frac{1}{\gamma \log^2 M} \prod_{p|\gamma} \left(1 + O\left(\frac{1}{p}\right)\right). \tag{6.3}$$

Notice that (6.2) and (6.3) ensure that $\lambda(l) \ll l^{-1+\varepsilon}$.

With these conventions in mind, we proceed to evaluate the first and second mollified moments.

6.1. *The first mollified moment.* Using Lemma 2.2 and Proposition 1.1 we see that

$$\mathcal{S}(M(d)L(\tfrac{1}{2}, \chi_{8d}); \Phi) = 2\mathcal{S}_M(M(d)A_1(d); \Phi) + O(X^{-\varepsilon}).$$

Using Proposition 1.2 we get

$$2\mathcal{S}_M(M(d)A_1(d); \Phi)$$
$$= \frac{C}{\zeta(2)} \hat{\Phi}(0) \sum_{l \leq M} \frac{\lambda(l)}{g(l)} \left(\log \frac{\sqrt{X}}{l} + C_2 + \sum_{p|l} \frac{C_2(p)}{p} \log p\right) + O(X^{-\varepsilon}).$$

Define $g_1(\gamma)$ to be the multiplicative function defined on primes by

$$g_1(p) = \frac{1}{g(p)} - \frac{2p}{h(p)(p+1)}.$$

It's easy to see that $g_1(p) = -1 + O(\frac{1}{p})$. Writing $\lambda$ in terms of $\xi$ using (6.2) we deduce that

$$\sum_l \frac{\lambda(l)}{g(l)} \log \frac{\sqrt{X}}{l} = \sum_l \sum_a \xi(la) \frac{1}{g(l)} \frac{\mu(a)ad(a)}{h(a)\sigma(a)} \log \frac{\sqrt{X}}{l}$$
$$= \sum_\gamma \xi(\gamma) \sum_{al=\gamma} \frac{1}{g(l)} \frac{\mu(a)ad(a)}{h(a)\sigma(a)} \log \frac{\sqrt{X}}{l}$$
$$= \sum_\gamma \xi(\gamma) g_1(\gamma) \left(\log(\sqrt{X}\gamma) + O\left(\sum_{p|\gamma} \frac{\log p}{p}\right)\right).$$

By (6.3) this is
$$\sum_\gamma \xi(\gamma) g_1(\gamma)(\log \sqrt{X}\gamma) + O\left(\frac{1}{\log X}\right).$$

Similarly one sees that
$$\sum_l \frac{\lambda(l)}{g(l)}\left(C_2 + \sum_{p|l} \frac{\log p}{p} C_2(p)\right) \ll \frac{1}{\log X}.$$

We have shown that the first mollified moment is

(6.4) $$\frac{C}{\zeta(2)}\hat{\Phi}(0) \sum_\gamma \xi(\gamma) g_1(\gamma) \log(\sqrt{X}\gamma) + O\left(\frac{1}{\log X}\right).$$

6.2. *The second mollified moment.* Using Lemma 2.2, Proposition 1.1 and Proposition 1.3 we know that the second mollified moment is (with an error $O(X^{-\varepsilon})$)

$$\frac{D\hat{\Phi}(0)}{36\zeta(2)} \sum_l \left(\sum_{rs=l} \lambda(r)\lambda(s)\right) \frac{\sqrt{l}}{h(l)} \frac{d(l_1)}{\sqrt{l_1}} \frac{l_1}{\sigma(l_1)} \left(\log^3\left(\frac{X}{l_1}\right) - 3\sum_{p|l_1} \log^2 p \log\left(\frac{X}{l_1}\right) + \mathcal{O}(l)\right).$$

We write $r = a\alpha$ and $s = b\alpha$ where $a$ and $b$ are coprime. Since we assumed that $\lambda$ is supported on square-frees we note that $\alpha = l_2$ and $l_1 = ab$. Thus the above may be rewritten as

$$\frac{D\hat{\Phi}(0)}{36\zeta(2)} \sum_\alpha \frac{\alpha}{h(\alpha)} \sum_{\substack{a,b \\ (a,b)=1}} \frac{\lambda(a\alpha)}{h(a)} \frac{\lambda(b\alpha)}{h(b)} \frac{ad(a)}{\sigma(a)} \frac{bd(b)}{\sigma(b)}$$
$$\times \left(\log^3\left(\frac{X}{ab}\right) - 3\sum_{p|ab} \log^2 p \log\left(\frac{X}{ab}\right) + \mathcal{O}(\alpha^2 ab)\right).$$

Since $\sum_{\beta|(a,b)} \mu(\beta) = 1$ or 0 depending on whether $(a,b)=1$ or not, the above becomes
(6.5)
$$\frac{D\hat{\Phi}(0)}{36\zeta(2)} \sum_\alpha \frac{\alpha}{h(\alpha)} \sum_\beta \frac{\mu(\beta)}{h(\beta)^2} \frac{\beta^2 d(\beta)^2}{\sigma(\beta)^2} \sum_{a,b} \frac{\lambda(a\alpha\beta)}{h(a)} \frac{\lambda(b\alpha\beta)}{h(b)} \frac{ad(a)}{\sigma(a)} \frac{bd(b)}{\sigma(b)}$$
$$\times \left(\log^3\left(\frac{X}{ab\beta^2}\right) - 3\sum_{p|ab\beta} \log^2 p \log\left(\frac{X}{ab\beta^2}\right) + \mathcal{O}(\alpha^2\beta^2 ab)\right).$$

We now define a multiplicative function $H(n)$ by setting
$$H(p) = 1 - \frac{4p}{h(p)(p+1)^2} = 1 + O\left(\frac{1}{p}\right).$$

Observe that, for nonnegative integers $j$, and square-free integers $\gamma$,

$$\sum_{\alpha\beta=\gamma} \frac{\alpha}{h(\alpha)} \frac{\mu(\beta)}{h(\beta)^2} \frac{\beta^2 d(\beta)^2}{\sigma(\beta)^2} (\log \beta)^j = \frac{\gamma}{h(\gamma)} \sum_{m|\gamma} \Lambda_j(m) \frac{\mu(m)}{h(m)} \frac{m d(m)^2}{\sigma(m)^2} H\left(\frac{\gamma}{m}\right)$$

$$\ll_j \frac{\gamma}{h(\gamma)} H(\gamma) \sum_{m|\gamma} \frac{\Lambda_j(m)}{m}.$$

Further, using (6.3), we note that for any square-free integer $\gamma$ and any nonnegative integer $j$,

$$\sum_a \frac{\lambda(a\gamma)}{h(a)} \frac{a d(a)}{\sigma(a)} (\log a)^j = \sum_m \Lambda_j(m) \sum_{\substack{a \\ m|a}} \frac{\lambda(a\gamma)}{h(a)} \frac{a d(a)}{\sigma(a)}$$

$$= \sum_m \Lambda_j(m) \frac{m d(m)}{h(m)\sigma(m)} \xi(m\gamma)$$

$$\ll_j \frac{1}{\gamma \log^2 M} \prod_{p|\gamma}\left(1 + O\left(\frac{1}{p}\right)\right) \sum_{m \leq M/\gamma} \frac{\Lambda_j(m)}{m}$$

$$\ll_j \frac{(\log M)^{j-2}}{\gamma} \prod_{p|\gamma}\left(1 + O\left(\frac{1}{p}\right)\right).$$

From these two observations we see easily that

$$\sum_\alpha \frac{\alpha}{h(\alpha)} \sum_\beta \frac{\mu(\beta)}{h(\beta)^2} \frac{\beta^2 d(\beta)^2}{\sigma(\beta)^2}$$

$$\times \sum_{a,b} \frac{\lambda(a\alpha\beta)}{h(a)} \frac{\lambda(b\alpha\beta)}{h(b)} \frac{a d(a)}{\sigma(a)} \frac{b d(b)}{\sigma(b)} \left(\log^3 \frac{X}{ab\beta^2} - \log^3 \frac{X}{ab}\right)$$

$$\ll \sum_{j=1}^3 \sum_{\gamma \leq M} \frac{\gamma H(\gamma)}{h(\gamma)} \sum_{m|\gamma} \frac{\Lambda_j(m)}{m} \frac{(\log X)^{-1-j}}{\gamma^2} \prod_{p|\gamma}\left(1 + O\left(\frac{1}{p}\right)\right) \ll \frac{1}{\log X}.$$

Similarly we see that

$$\sum_\alpha \frac{\alpha}{h(\alpha)} \sum_\beta \frac{\mu(\beta)}{h(\beta)^2} \frac{\beta^2 d(\beta)^2}{\sigma(\beta)^2} \sum_{a,b} \frac{\lambda(a\alpha\beta)}{h(a)} \frac{\lambda(b\alpha\beta)}{h(b)} \frac{a d(a)}{\sigma(a)} \frac{b d(b)}{\sigma(b)}$$

$$\times \left(\log \frac{X}{ab}\left(\sum_{p|a} \log^2 p + \sum_{p|b} \log^2 p\right) - \log \frac{X}{ab\beta^2} \sum_{p|ab\beta} \log^2 p\right) \ll \frac{1}{\log X},$$

and that

$$\sum_\alpha \frac{\alpha}{h(\alpha)} \sum_\beta \frac{\mu(\beta)}{h(\beta)^2} \frac{\beta^2 d(\beta)^2}{\sigma(\beta)^2}$$

$$\times \sum_{a,b} \frac{\lambda(a\alpha\beta)}{h(a)} \frac{\lambda(b\alpha\beta)}{h(b)} \frac{a d(a)}{\sigma(a)} \frac{b d(b)}{\sigma(b)} \mathcal{O}(\alpha^2 \beta^2 ab) \ll \frac{1}{\log X}.$$

We have shown that the second mollified moment is (with an error $O(1/\log X)$)

(6.6)
$$\frac{D\hat{\Phi}(0)}{36\zeta(2)} \sum_\gamma \frac{\gamma H(\gamma)}{h(\gamma)} \sum_{a,b} \frac{\lambda(a\gamma)}{h(a)} \frac{\lambda(b\gamma)}{h(b)} \frac{ad(a)}{\sigma(a)} \frac{bd(b)}{\sigma(b)} \left( \log^3 \frac{X}{ab} - 3\log \frac{X}{ab} \left( \sum_{p|a} \log^2 p + \sum_{p|b} \log^2 p \right) \right).$$

### 6.3. Completion of the proof.
Roughly speaking our expression for the second mollified moment looks like

(6.7)
$$\frac{D\hat{\Phi}(0)}{36\zeta(2)} \log^3 X \sum_\gamma \frac{\gamma H(\gamma)}{h(\gamma)} \xi(\gamma)^2.$$

This is a diagonal quadratic form in the $\xi$'s and we shall choose our mollifier so as to minimize (6.7) for fixed (6.4). Obviously this is achieved by choosing $\xi(\gamma)$ (for odd square-free $\gamma \leq M$) to be proportional to

$$\frac{h(\gamma)g_1(\gamma)}{\gamma H(\gamma)} \log(\sqrt{X}\gamma).$$

In fact, we shall choose (for odd square-free $\gamma \leq M$)

(6.8)
$$\xi(\gamma) = \frac{C}{D \log^3 M} \frac{h(\gamma)g_1(\gamma)}{\gamma H(\gamma)} \log(\sqrt{X}\gamma).$$

Observe that our choice (6.8) meets the constraint (6.3) imposed earlier.

An elementary argument shows that

(6.9)
$$\frac{C^2}{D} \sum_{\gamma \leq x} \mu^2(2\gamma) \frac{h(\gamma)g_1(\gamma)^2}{\gamma H(\gamma)} = \frac{C^2}{D} \frac{1}{2} \prod_{p \geq 3} \left(1 - \frac{1}{p}\right) \left(1 + \frac{h(p)g_1(p)^2}{pH(p)}\right) (\log x + O(1))$$
$$= \frac{4}{9}(\log x + O(1)).$$

From this and partial summation we get that the first mollified moment is

(6.10)
$$\sim \frac{C^2}{D\zeta(2)} \frac{\hat{\Phi}(0)}{\log^3 M} \sum_{\gamma \leq M} \mu(2\gamma)^2 \frac{h(\gamma)g_1(\gamma)^2}{\gamma H(\gamma)} \log^2(\sqrt{X}\gamma)$$
$$\sim \frac{2}{9} \left(\left(1 + \frac{1}{\theta}\right)^3 - \frac{1}{\theta^3}\right) \frac{2\hat{\Phi}(0)}{3\zeta(2)},$$

For nonnegative integers $j$ note that

$$\xi_j(\gamma) := \sum_a \frac{\lambda(a\gamma)}{h(a)} \frac{ad(a)}{\sigma(a)} (\log a)^j = \sum_a \frac{\lambda(a\gamma)}{h(a)} \frac{ad(a)}{\sigma(a)} \sum_{m|a} \Lambda_j(m)$$
$$= \sum_m \Lambda_j(m) \frac{md(m)}{h(m)\sigma(m)} \xi(m\gamma).$$

Note that $\xi_j(\gamma)$ is supported only on odd square-free integers $\leq M$, and that for such a $\gamma$ our choice (6.8) gives

$$\xi_j(\gamma) = \frac{C}{D\log^3 M} \frac{h(\gamma)g_1(\gamma)}{\gamma H(\gamma)} \sum_{\substack{m \leq M/\gamma \\ (m,\gamma)=1}} \mu(m)d(m)\frac{\Lambda_j(m)}{m}$$

$$\times \prod_{p|m}\left(1+O\left(\frac{1}{p}\right)\right)\log(\sqrt{X}m\gamma).$$

It is easy to check that

$$\sum_{\substack{m \leq x \\ (m,\gamma)=1}} \mu(m)d(m)\frac{\Lambda_j(m)}{m}\prod_{p|m}\left(1+O\left(\frac{1}{p}\right)\right)$$

$$= \sum_{m \leq x} \mu(m)d(m)\frac{\Lambda_j(m)}{m}\prod_{p|m}\left(1+O\left(\frac{1}{p}\right)\right) + O\left(\sum_{q|\gamma} 2^j \sum_{\substack{m \leq x \\ q|m}} \frac{\Lambda_j(m)}{m}\prod_{p|m}\left(1+O\left(\frac{1}{p}\right)\right)\right)$$

$$= \begin{cases} -2\log x + O(1 + \sum_{q|\gamma} \frac{\log q}{q}) & \text{if } j = 1, \\ \log^2 x + O(\log x(1 + \sum_{q|\gamma} \frac{\log q}{q})) & \text{if } j = 2, \\ \ll \log^2 x (1 + \sum_{q|\gamma} \frac{\log q}{q}) & \text{if } j = 3. \end{cases}$$

From this and partial summation we get

(6.11a)
$$\xi_1(\gamma) = -\frac{C}{D\log^3 M}\frac{h(\gamma)g_1(\gamma)}{\gamma H(\gamma)}\left(2\log\left(\frac{M}{\gamma}\right)\log(\sqrt{X}\gamma) + \log^2\left(\frac{M}{\gamma}\right)\right.$$
$$\left. + O\left(\log M\left(1+\sum_{q|\gamma}\frac{\log q}{q}\right)\right)\right),$$

(6.11b)
$$\xi_2(\gamma) = \frac{C}{D\log^3 M}\frac{h(\gamma)g_1(\gamma)}{\gamma H(\gamma)}\left(\log^2\left(\frac{M}{\gamma}\right)\log(\sqrt{X}\gamma) + \frac{2}{3}\log^3\left(\frac{M}{\gamma}\right)\right.$$
$$\left. + O\left(\log^2 M\left(1+\sum_{q|\gamma}\frac{\log q}{q}\right)\right)\right),$$

and

(6.11c)
$$\xi_3(\gamma) \ll \frac{|h(\gamma)g_1(\gamma)|}{\gamma H(\gamma)}\left(1+\sum_{q|\gamma}\frac{\log q}{q}\right).$$

Expanding $\log^3(X/ab)$ in terms of $\log X$, $\log a$ and $\log b$, we may write

(6.12) $$\frac{D\hat{\Phi}(0)}{36\zeta(2)} \sum_\gamma \frac{\gamma H(\gamma)}{h(\gamma)} \sum_{a,b} \frac{\lambda(a\gamma)}{h(a)}\frac{ad(a)}{\sigma(a)}\frac{\lambda(b\gamma)}{h(b)}\frac{bd(b)}{\sigma(b)}\log^3\frac{X}{ab}$$

as a linear combination of terms
$$\frac{D\hat{\Phi}(0)}{36\zeta(2)} \sum_\gamma \frac{\gamma H(\gamma)}{h(\gamma)} \xi_j(\gamma)\xi_k(\gamma)(\log X)^l, \qquad \text{where } j+k+l=3.$$

These terms may be evaluated by appealing to (6.11a,b,c) and then using (6.9) and partial summation. In this manner we show that

$$(6.12) \sim \left(\frac{2}{81} + \frac{28}{135\theta} + \frac{11}{18\theta^2} + \frac{70}{81\theta^3} + \frac{16}{27\theta^4} + \frac{4}{27\theta^5}\right)\frac{2\hat{\Phi}(0)}{3\zeta(2)}.$$

This handles one of the terms in our asymptotic formula (6.6) for the second mollified moment.

To handle the other term, we note that for odd, square-free $\gamma \leq M$

$$\sum_a \frac{\lambda(a\gamma)}{h(a)} \frac{ad(a)}{\sigma(a)} \sum_{p|a} \log^2 p$$

$$= \sum_p \log^2 p(2 + O(p^{-1}))\xi(\gamma p)$$

$$= -\frac{C}{D}\frac{h(\gamma)g_1(\gamma)}{\gamma H(\gamma)\log^3 M} \sum_{p \leq M/\gamma} \frac{2\log^2 p}{p}\left(1 + O\left(\frac{1}{p}\right)\right)\log(\sqrt{X}\gamma p)$$

$$= -\frac{C}{D}\frac{h(\gamma)g_1(\gamma)}{\gamma H(\gamma)\log^3 M}\left(\log^2 \frac{M}{\gamma}\log(\sqrt{X}\gamma) + \frac{2}{3}\log^3 \frac{M}{\gamma} + O(\log^2 X)\right),$$

and similarly

$$\sum_a \frac{\lambda(a\gamma)}{h(a)} \frac{ad(a)}{\sigma(a)} \log a \sum_{p|a} \log^2 p$$

$$= \sum_p \log^2 p \sum_{p|a} \frac{\lambda(a\gamma)}{h(a)} \frac{ad(a)}{\sigma(a)} \log a$$

$$= \sum_p \log^2 p(2 + O(p^{-1})) \sum_a \frac{\lambda(ap\gamma)}{h(a)} \frac{ad(a)}{\sigma(a)} \log(ap)$$

$$= \sum_p 2\log^2 p \left(1 + O\left(\frac{1}{p}\right)\right)\left(\log p\, \xi(p\gamma) + \xi_1(p\gamma)\right).$$

Here we employ (6.8) and (6.11a), and use partial summation. It transpires that the main terms cancel out, and we find that the above is

$$\ll \frac{|h(\gamma)g_1(\gamma)|}{\gamma H(\gamma)}\left(1 + \sum_{q|\gamma} \frac{\log q}{q}\right).$$

Using these we find that

$$-\frac{D\hat{\Phi}(0)}{12\zeta(2)}\sum_{\gamma}\frac{\gamma H(\gamma)}{h(\gamma)}\sum_{a,b}\frac{\lambda(a\gamma)}{h(a)}\frac{ad(a)}{\sigma(a)}\frac{\lambda(b\gamma)}{h(b)}\frac{bd(b)}{\sigma(b)}\log\frac{X}{ab}\left(\sum_{p|a}\log^2 p+\sum_{p|b}\log^2 p\right)$$

$$\sim \left(\frac{2}{81}+\frac{4}{45\theta}+\frac{7}{54\theta^2}+\frac{2}{27\theta^3}\right)\frac{2\hat{\Phi}(0)}{3\zeta(2)}.$$

Combining this with our evaluation of (6.12) we find that the second mollified moment is

$$(6.13)\qquad \sim \left(\frac{4}{81}+\frac{8}{27\theta}+\frac{20}{27\theta^2}+\frac{76}{81\theta^3}+\frac{16}{27\theta^4}+\frac{4}{27\theta^5}\right)\frac{2\hat{\Phi}(0)}{3\zeta(2)}.$$

We choose $\Phi$ to be an approximation to the characteristic function of $(1,2)$ so that $\hat{\Phi}(0) \sim 1$. By Cauchy's inequality, and the evaluations of the mollified moments,

$$\sum_{\substack{X\leq d\leq 2X \\ d \text{ odd} \\ L(\frac{1}{2},\chi_{8d})\neq 0}} \mu(d)^2 \geq \sum_{\substack{d \text{ odd} \\ L(\frac{1}{2},\chi_{8d})\neq 0}} \mu(d)^2 \Phi\left(\frac{d}{X}\right) \geq X\frac{\mathcal{S}(L(\frac{1}{2},\chi_{8d})M(d);\Phi)^2}{\mathcal{S}(L(\frac{1}{2},\chi_{8d})^2 M(d)^2;\Phi)}$$

$$\geq \left(1-\frac{1}{(\theta+1)^3}\right)\frac{4}{\pi^2}X = \left(\frac{7}{8}+o(1)\right)\sum_{\substack{X\leq d\leq 2X \\ d \text{ odd}}} \mu(d)^2,$$

upon taking $\theta = 1 - \varepsilon$. Take this with $X = x/2$, $x/4$, …, and sum to get Theorem 1.

## 7. Sketch proof of Theorem 2

By Lemma 2.2 and by modifying the proof of Proposition 1.1, we have

$$\mathcal{S}(L(\tfrac{1}{2},\chi_{8d})^3;\Phi) = 2\mathcal{S}_M(A_3(d);\Phi) + O\left(\frac{X^\varepsilon}{Y}\right)$$

$$= 2\sum_{n=1}^{\infty}\frac{d_3(n)}{\sqrt{n}}\mathcal{S}_M\left(\left(\frac{8d}{n}\right);H_n\right) + O\left(\frac{X^\varepsilon}{Y}\right),$$

where $H_n(t) = \Phi(t)\omega_3(n\pi^{\frac{3}{2}}/(8Xt)^{\frac{3}{2}})$. By Poisson summation (Lemma 2.6 above) this becomes

$$\sum_{\substack{n=1 \\ (n,2)=1}}^{\infty}\frac{d_3(n)}{n^{\frac{3}{2}}}\sum_{\substack{\alpha\leq Y \\ (\alpha,2n)=1}}\frac{\mu(\alpha)}{\alpha^2}\sum_{k=-\infty}^{\infty}(-1)^k G_k(n)\tilde{H}_n\left(\frac{kX}{2\alpha^2 n}\right).$$

The method of Section 5.1 shows that the $k = 0$ term above contributes

$$Q_{0,\Phi}(\log X) + O(X^\varepsilon Y^{-1}),$$

for a polynomial $Q_0$ of degree 6 whose coefficients involve linear combinations of the parameters $\check{\Phi}^{(j)}(0)$ for $j = 0, \ldots, 6$. As in Section 5.2 we may extract a secondary principal contribution from the $k = \square$ terms. These may be evaluated as in Section 5.3 to give a contribution

$$Q_{1,\Phi}(\log X) + O(X^\varepsilon Y^{-1}),$$

for a polynomial $Q_1$ also of degree 6 whose coefficients are again linear combinations of the $\check{\Phi}^{(j)}(0)$ for $j = 0, \ldots, 6$. Lastly the remainder terms (arising from $k \neq 0, \square$ values) are estimated analogously to Section 5.4, and are $\ll \Phi_{(2)} \Phi_{(3)}^\varepsilon Y X^{-\frac{1}{4}+\varepsilon}$. Thus one may show that

$$\mathcal{S}(L(\tfrac{1}{2}, \chi_{8d})^3; \Phi) = Q_{2,\Phi}(\log X) + O\left(\frac{X^\varepsilon}{Y} + \Phi_{(2)} \Phi_{(3)}^\varepsilon Y X^{-\frac{1}{4}}\right),$$

for a polynomial $Q_{2,\Phi} = Q_{1,\Phi} + Q_{0;\Phi}$ of degree 6.

We now argue exactly as in Section 5.5, choosing $\Phi$ to be a good approximation to the characteristic function of $(1,2)$ (with $Y = Z = X^{\frac{1}{12}}$), thus proving Theorem 2.

INSTITUTE FOR ADVANCED STUDY, PRINCETON NJ
*E-mail address*: ksound@ias.edu